\documentclass[10pt]{amsart}

\usepackage{amsfonts,amssymb,amsthm,amsmath}
\usepackage{eucal}
\usepackage{ifthen,epsfig,color,graphics,rotating,hangcaption}

\newtheorem{thm}{Theorem}[section]
\newtheorem{prop}[thm]{Proposition}
\newtheorem{cor}[thm]{Corollary}
\newtheorem{lem}[thm]{Lemma}

\theoremstyle{definition}
\newtheorem{defn}[thm]{Definition}
\newtheorem{exmp}[thm]{Example}
\newtheorem{alg}[thm]{Algorithm}

\theoremstyle{remark}

\newtheorem{rem}{Remark}

\renewcommand{\(}{\left(}
\renewcommand{\)}{\right)}  
\newcommand{\abs}[1]{\left\arrowvert #1\right\arrowvert}

\newcommand{\onorm}[1]{\abs{#1}}

\newcommand{\sonorm}[1]{\onorm{#1}}

\newcommand{\eonorm}[1]{\onorm{#1}_{e}}

\newcommand{\ponorm}[1]{\onorm{#1}_{\rho}}
\newcommand{\konorm}[1]{\onorm{#1}_{k}}
\newcommand{\jonorm}[1]{\onorm{#1}_{j}}
\newcommand{\Phionorm}[1]{\onorm{#1}_{\Phi}}

\newcommand{\nbd}[2]{\mathcal{N}(#1, #2)}

\newcommand{\field}[1]{\mathbb{#1}}
\newcommand{\DD}{\ensuremath{\field{D}}} 
\newcommand{\FF}{\ensuremath{\field{F}}}
\newcommand{\KK}{\ensuremath{\field{K}}} 
\newcommand{\II}{\ensuremath{\field{I}}} 
 
\newcommand{\CC}{\ensuremath{\field{C}}} 
\newcommand{\Ct}{\ensuremath{\field{C}^2}}

\newcommand{\RR}{\ensuremath{\field{R}}}
\newcommand{\Rt}{\ensuremath{\field{R}^2}}
\newcommand{\Rn}{\ensuremath{\field{R}^n}}

\newcommand{\boxcov}{\mathcal{B}}

\newcommand{\boxchcn}{box chain construction}

\newcommand{\boxmod}{box chain model}
\newcommand{\boxchmod}{box chain model}
\newcommand{\boxchmods}{box chain models}

\newcommand{\inKselsub}{sink basin subdivision}
\newcommand{\wkcycsub}{weak cycle subdivision}

\newcommand{\boxmet}{box metric}
\newcommand{\namealg}{Handicap Hedging Algorithm}
\newcommand{\boxexp}{box expansive}
\newcommand{\boxJ}{box Julia set}

\newcommand{\drawfighyponeF}{\scalebox{.35}{\includegraphics{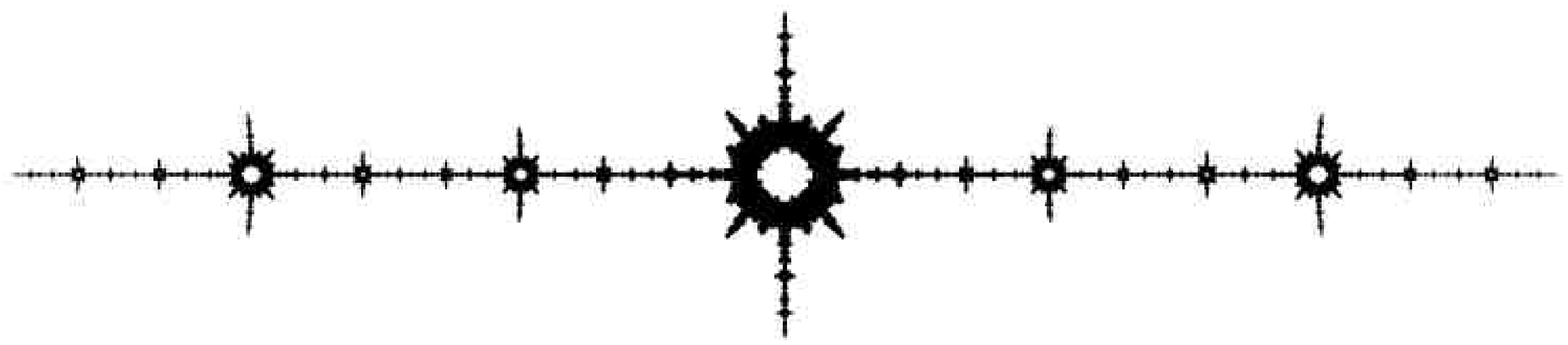}}}
\newcommand{\drawfigcantcaulFA}{\scalebox{.225}{\includegraphics{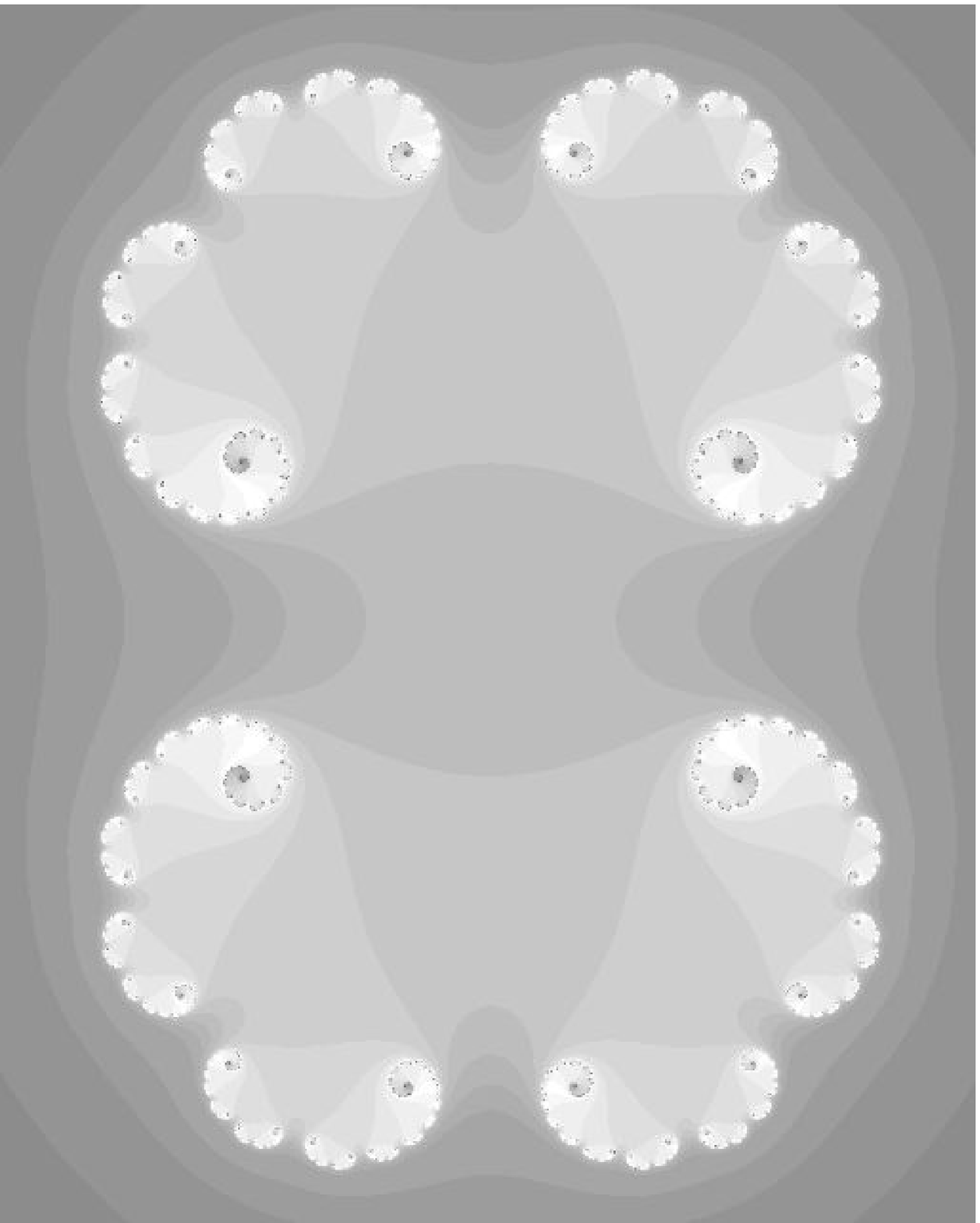}}}
\newcommand{\drawfigcantcaulmet}{\scalebox{.35}{\includegraphics{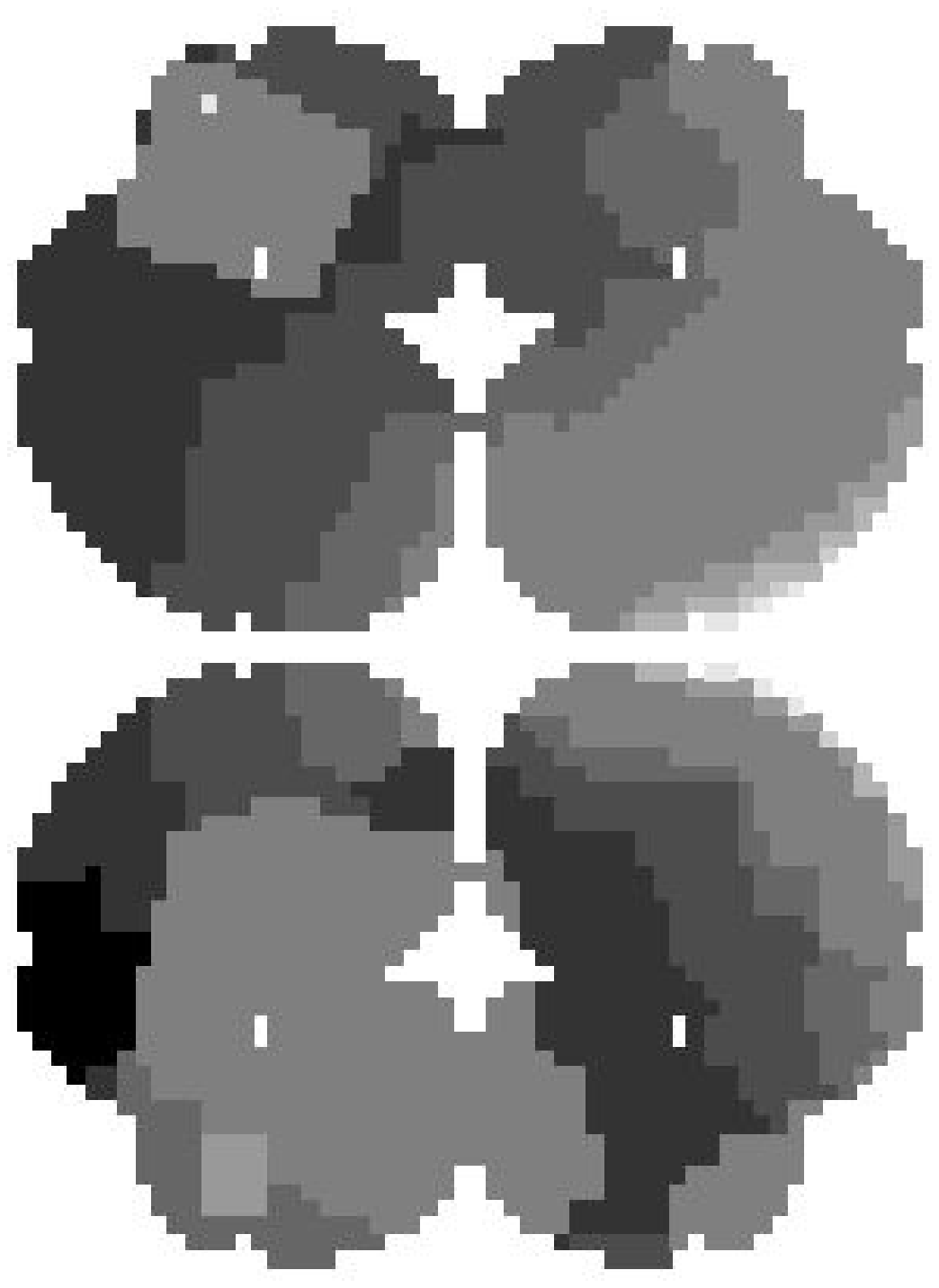}}}

\newcommand{\drawfigbasFA}{\scalebox{.325}{\includegraphics{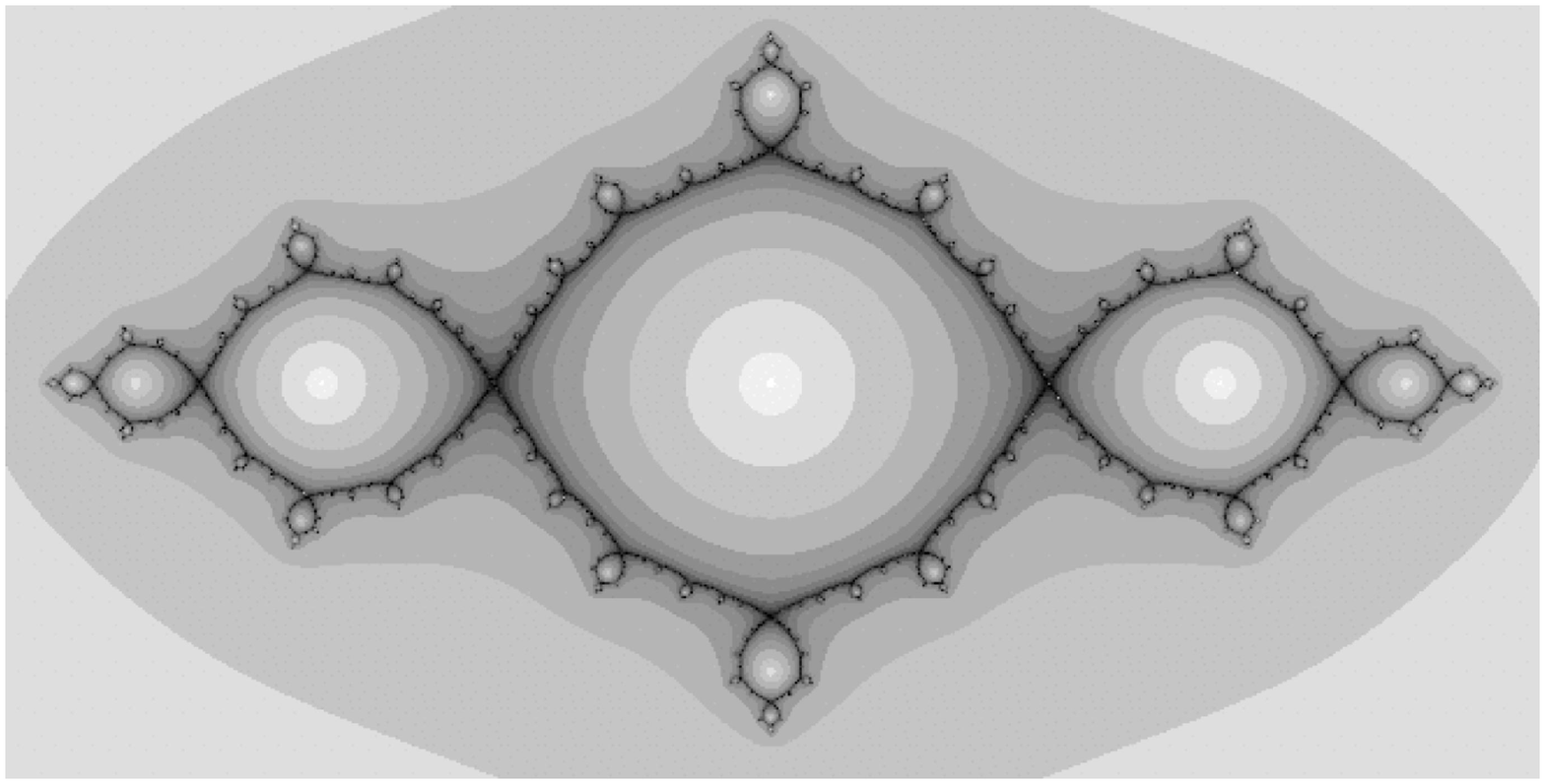}}}

\newcommand{\drawfigfirstmetex}{\scalebox{2}{\includegraphics{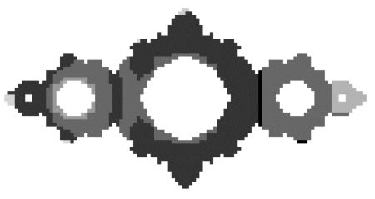}}}
\newcommand{\drawfigcubrabtwo}{\scalebox{.4}{\includegraphics{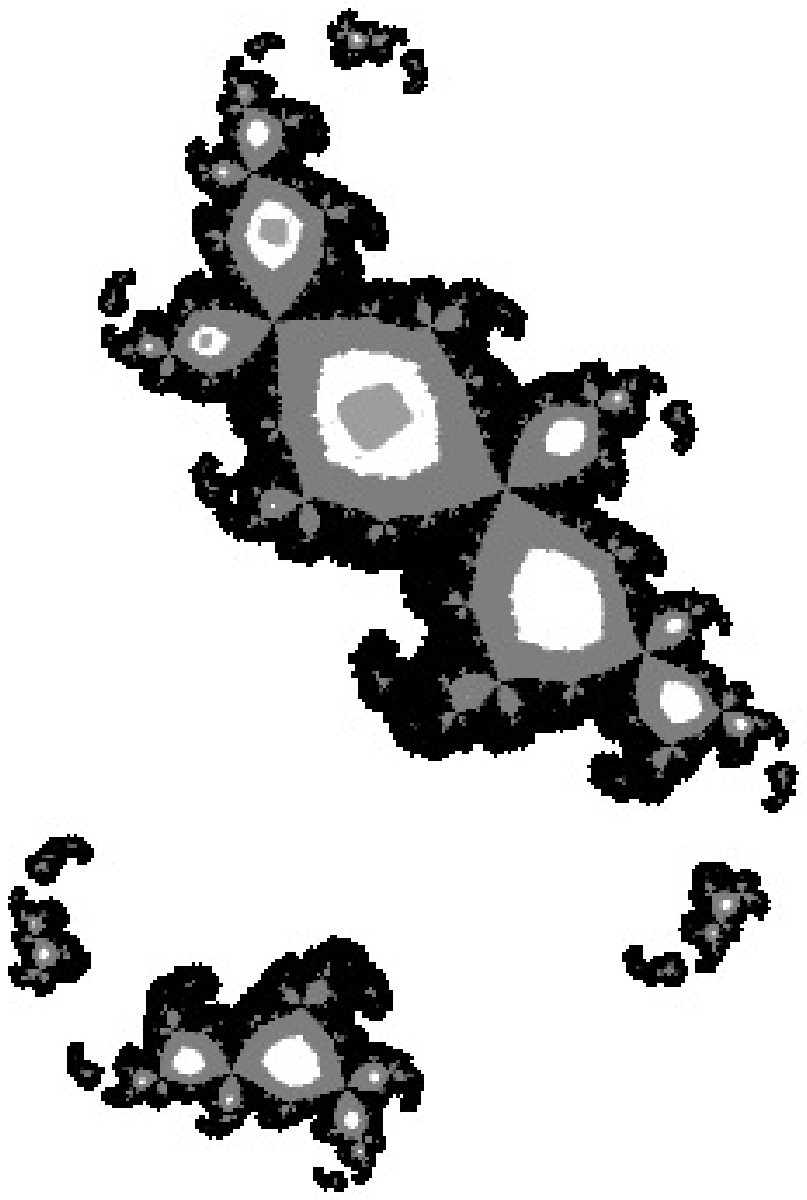}}}
\newcommand{\drawfigcubrabthreea}{\scalebox{.425}{\includegraphics{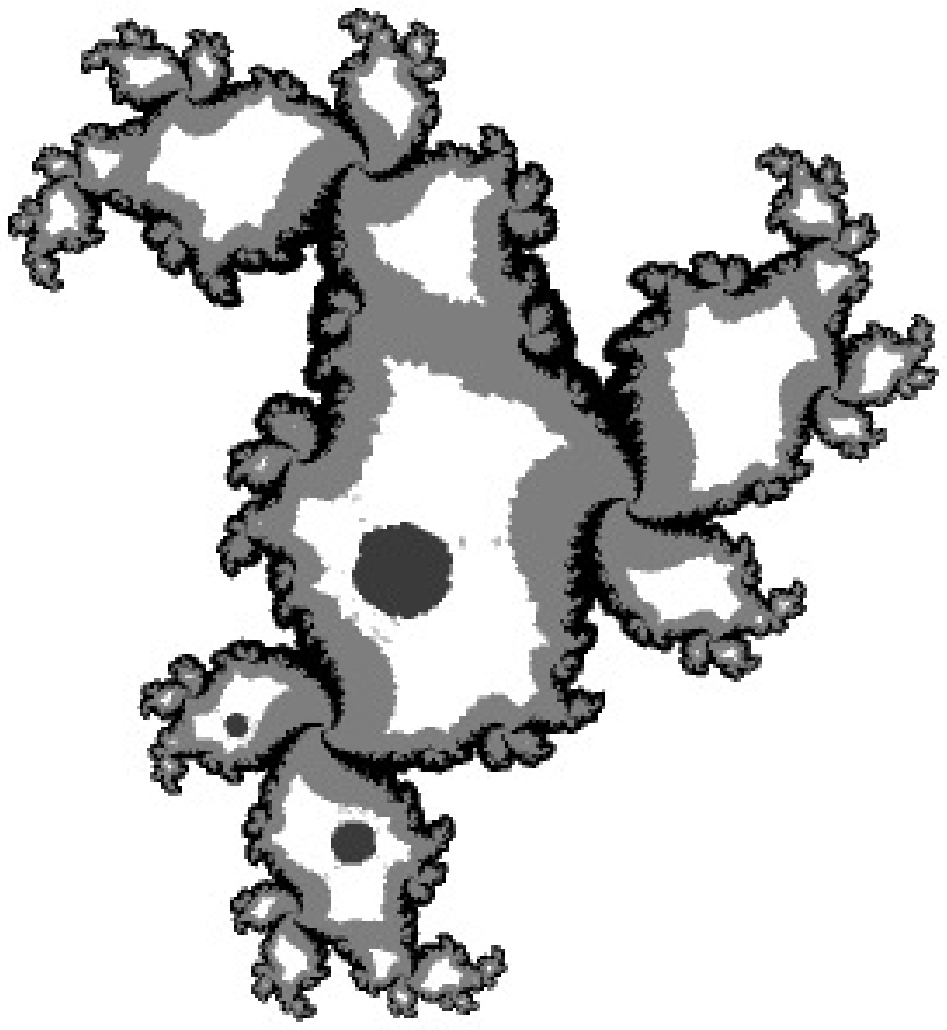}}}
\newcommand{\drawfigcubrabthreeb}{\scalebox{.425}{\includegraphics{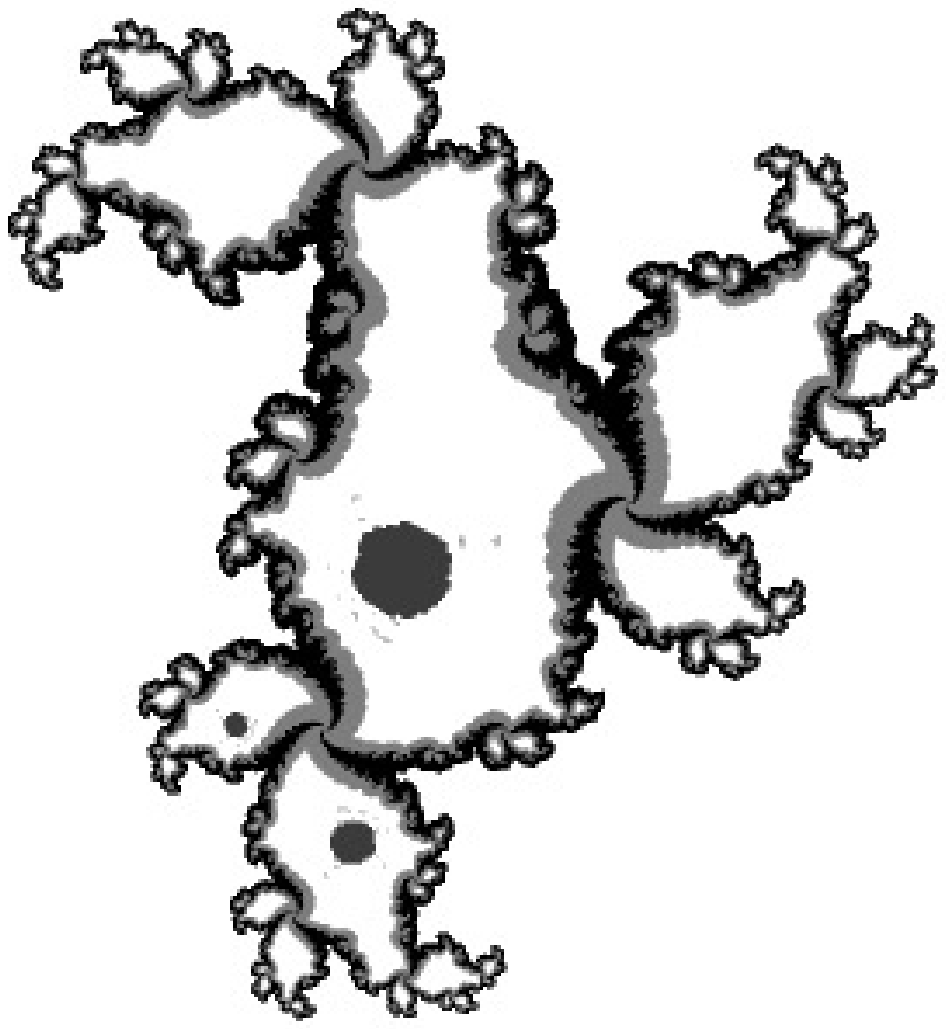}}}

\newcommand{\drawfiggraphzero}{\includegraphics{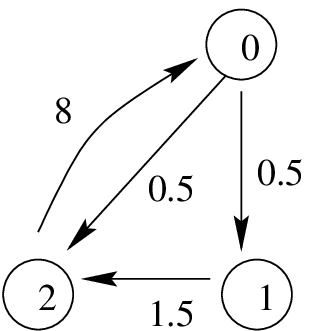}}

\begin{document}

\title[Constructing an expanding metric in one complex dimension]{Constructing an expanding metric\\ for dynamical systems in one complex variable}
\author[S. L. ~Hruska]{Suzanne Lynch Hruska}
\address{Department of Mathematics\\
Indiana University\\
Rawles Hall\\
Bloomington, IN 47405, USA
}
\email{shruska@msm.umr.edu}

\date{\today}

\begin{abstract}    
We describe a rigorous computer algorithm for attempting to construct an explicit, discretized metric for which a   polynomial map $f \colon \CC \rightarrow \CC$ is expansive on a neighborhood of the Julia set, $J$.  We show    construction of such a metric proves the map is hyperbolic.  We also examine the question of whether the algorithm can be improved, and the related question of how to build a metric close to euclidean.  Finally, we give several examples generated with our implementation of this algorithm.
\end{abstract}


\renewcommand{\subjclassname}{\textup{2000} Mathematics Subject Classification}


\maketitle

\section{Introduction}
\label{sec:intro}

Our main concern in this paper is 
to develop and use \textit{rigorous} computer investigations to study the dynamics of polynomial maps of $\CC$ of degree $d>1$; for example, $f(z) = z^2 + c$. 
For complex polynomial maps, the invariant set of interest is the {\em Julia set}, $J$.  The Julia set can be defined as the topological boundary of the set of points with bounded orbits.  
Intuitively, $J$ is precisely where the chaotic dynamics occurs.  For example, $f$ is topologically transitive on $J$; also, $J$ is non-empty and perfect (see for instance \cite{Mil}).
 
The hyperbolic polynomials are a large class of maps with chaotic dynamics, but whose
stability properties make them amenable to computer study.
A polynomial map $f$ is called \textit{hyperbolic}, or \textit{expansive}, if $f$ is uniformly expanding on some neighborhood of $J$, with respect to some riemannian metric.   Uniform expansion forces some
dynamical rigidity.
For example, 
hyperbolicity of a polynomial map $f$ implies \textit{structural stability}, \textit{i.e.}, in a neighborhood of $f$ in parameter space, the dynamical behavior is of constant topological conjugacy type.  
Further,  if $f$ is hyperbolic then the orbit of every point in the complement of $J$ tends to either some attracting periodic orbit or infinity.  
In fact, a polynomial map $f$ is hyperbolic if and only if the orbit of every  {\em critical} point of $f$ tends to an attracting periodic orbit or infinity.  
Thus the fate of the critical points provides a straightforward test for hyperbolicity.  In this paper we develop an alternate test for hyperbolicity which produces  more explicit information about the dynamics of any given map. 


We begin with the work of \cite{SLHone} as a foundation.  There we described a 
rigorous algorithm (and its implementation) for constructing a neighborhood $\boxcov$ of 
$J$, 
and a graph $\Gamma$ which models the dynamics of $f$ on $\boxcov$.
In this paper, we further our program by developing a rigorous algorithm for attempting to construct a metric in which a given $f$ is expansive, by some
uniform factor $L>1$, on the neighborhood $\boxcov$. We
show that successful construction of such a metric proves that $f$ is hyperbolic. 

In addition, we analyze 
both some limitations of and possible improvements to our basic algorithm.
We show, via counterexamples, that a significantly simpler algorithm will not suffice.  On the other hand, we do present an enhancement to our algorithm, to build a metric closer to euclidean.

One motivation for improving the algorithm and the metric which it builds is that in \cite{SLHthree}, we use this one dimensional algorithm 
as part of a computer-assisted proof of hyperbolicity of polynomial diffeomorphisms of $\Ct$, H\'{e}non mappings in particular.  The one dimensional algorithm is used in such a way that that the closer to euclidean we can build the metric, the more likely we are to succeed at proving hyperbolicity in $\Ct$.  Thus improving the current methods could lead to more interesting examples of proven hyperbolic polynomial diffeomorphisms of $\Ct$.

We implemented the algorithms of this paper in a program \textit{Hypatia}, 
\begin{footnote}{To obtain a copy of this unix program, write to the author.}\end{footnote}
and at the conclusion of the paper, we give several examples of metrics constructed which establish 
hyperbolicity of some polynomial maps of degrees two and three.

\subsection{Statement of main results}

In \cite{SLHone} we described the \textit{box chain construction}, for building a directed graph $\Gamma$ representing $f$ on some neighborhood $\boxcov$ of 
$J$. A similar approach is in the body of work described in the survey {\cite{KMisch} (cf \cite{Dell1,  Eiden, Osi, OsiCamp}).  
In this paper we need only:

\begin{defn} \label{defn:boxmodel}
Let $\Gamma = (\mathcal{V}, \mathcal{E})$ be a directed graph,
with vertex set $\mathcal{V} = \mathcal{V}(\Gamma) = \{ B_k \}_{k=1}^N$,  a finite collection of closed boxes in $\CC$, having disjoint interiors, and  such that the union of the boxes $\boxcov = \boxcov(\Gamma) = \cup_{k=1}^N B_k$ contains $J$.  Suppose there is a $\delta>0$ such that 
$\Gamma$ contains an edge from $B_k$ to $B_j$ if the image $f(B_k)$ intersects a $\delta$-neighborhood of $B_j$, \textit{i.e.},  
\[
\mathcal{E} \supset \{ (k, j) \colon f(B_k) \cap \mathcal{N} (B_j,{\delta} ) \neq \emptyset \}.
\]
Further, assume $\Gamma$ is strongly connected, \textit{i.e.}, for each pair of vertices $B_k, B_j$, there is a path in $\Gamma$ from $B_k$ to $B_j$, and vice-versa.   Then we call $\Gamma$ a \textit{{\boxchmod } of $f$ on $J$} and $\boxcov$ a \textit{\boxJ}.
\end{defn}


For the remainder of this paper, $\Gamma$ will denote a {\boxchmod } of $f$ on $J$, for some polynomial map $f \colon \CC \rightarrow \CC$ of degree greater than one.

\begin{rem} 
A {\boxmod } of $f$ on $J$ satisfies the definition of a {\em symbolic image of $f$}, given by Osipenko in \cite{Osiold}.
\end{rem}

Given a {\boxmod } $\Gamma$, 
our first goal in this paper is to define a version of hyperbolicity for $\Gamma$ which is checkable by computer, and which implies hyperbolicity for $f$.  
Our discrete version of hyperbolicity is parallel to the following standard
definition, from  \cite{Mil}:
a polynomial $f:\CC \to \CC$ is called \textit{hyperbolic} if there exists a
riemannian metric
$\mu$, on some
neighborhood $\boxcov$ of $J$, 
and an {\em
expansion constant} $\lambda > 1$,
  such that the derivative $D_zf$ at every point $z$ in $\boxcov$ satisfies
$ \onorm{D_zf(\mathbf{v})}_{\mu} \geq \lambda  \onorm{\mathbf{v}}_{\mu} $, for every
 vector $\mathbf{v}$ in the tangent space $T_z{\CC}$.

First, we introduce the following discrete version of a riemannian metric.

\begin{defn}
Let $\mathcal{V}= \{ B_k \}_{k=1}^N$ be a finite collection of closed boxes in $\CC$, 
with disjoint interiors.
Let $\Phi \colon (\sqcup_{k=1}^N B_k) \to \RR^+$ be a box constant function, \textit{i.e.}, 
for some set of positive constants $\{ \varphi_k \}_{k=1}^N$, we have
$\Phi(z) = \varphi_k$ if $z\in B_k$.
Then $\Phi$ times euclidean induces a metric on $T{\mathcal{V}} = \sqcup_{k=1}^N T B_k$.  
That is, for any $B_k \in \mathcal{V}$, $z \in B_k$, and $\mathbf{v} \in T_z B_k$, let $\konorm{\mathbf{v}} := \varphi_k \onorm{\mathbf{v}}$.  Let $\Phionorm{\cdot}$ be the induced norm on 
$\sqcup_{k=1}^N T B_k$. 
We call $\Phi$ a {\em \boxmet}, 
and say $\varphi_k$ is a {\em handicap} for $B_k$. %
\end{defn}

Note that if $z \in B_k \cap B_j$, then 
$\Phionorm{\mathbf{v}}=\konorm{\mathbf{v}}$ if we are considering $\mathbf{v} \in T_z B_k$,
but  $\Phionorm{\mathbf{v}}=\jonorm{\mathbf{v}}$ if we are considering $\mathbf{v} \in T_z B_j$.
We justify the choice of the name ``handicap'' a few paragraphs below.

Now we present our discretized version of hyperbolicity.

\begin{defn} 
\label{defn:boxexp}
Call $\Gamma  = (\mathcal{V}, \mathcal{E})$ {\em \boxexp }
if there exists a {\boxmet } $\Phi$ on $\mathcal{V}$ and a \textit{box expansion constant} $L > 1$ such that
 for all $(k,j)\in \mathcal{E}$, $z \in~B_k$, and $\mathbf{v} \in~T_zB_k$, we have 
$\jonorm{D_zf(\mathbf{v})} \geq L \konorm{\mathbf{v}}$. 
\end{defn}

The following theorem is a first piece of evidence that this definition is useful.

\begin{thm}  \label{thm:boxmetexp} 
Suppose there exists a {\boxmet } $\Phi$ and an $L>1$ for which $\Gamma$ is {\boxexp}.
Then $f$ is hyperbolic.

In particular, there exists a smooth function $\rho \colon \boxcov(\Gamma) \to \RR^+$, which defines
a riemannian metric on~$T {\boxcov}$  (by $\rho$ times euclidean), such that
for all $z \in \boxcov$, and $\mathbf{v} \in T_z \boxcov$, we have
$
\ponorm{D_zf(\mathbf{v})} \geq L \ponorm{\mathbf{v}}.
$
\end{thm}

The key step in the proof of this theorem is to smooth out the {\boxmet } $\Phi$ in a small neighborhood of the box boundaries, using a partition of unity argument and the ``edge overlap'' factor  $\delta$ from Definition~\ref{defn:boxmodel}.  The proof is given in Section~\ref{sec:met1exp}.

We develop an algorithm for attempting to build a {\boxmet } for which a given $\Gamma$ is {\boxexp}, called the \textit{\namealg}, in Section~\ref{sec:basicalg}.  Then we describe how any outcome of the algorithm gives useful dynamical information, in Section~\ref{sec:charboxexp}.  For example, we obtain:

\begin{thm} \label{thm:alg}
Given $\Gamma$ 
and $L>0$, the {\namealg } either
\begin{enumerate}
\item constructs a {\boxmet } for which $\Gamma$ is {\boxexp} by $L$, or 
\item produces a cycle of $\Gamma$ which is an obstruction, 
showing there exists no {\boxmet } for which $\Gamma$ is {\boxexp} by $L$.
\end{enumerate}
\end{thm}

Note that showing box expansion by $L=1$ may be 
instructive, but it is not enough to prove hyperbolicity.
Thus, if any algorithm, for example, the {\namealg}, builds a {\boxmet } for which some $\Gamma$ is
{\boxexp } by some $L>1$, then $f$ is proven hyperbolic. 

\subsection{Secondary results and discussion of the approach}

A natural question one might ask is whether there is a simple algorithm for showing box expansion, without constructing a {\boxmet};  for example, by examine the multipliers along the simple cycles in the graph.
(A cycle in a graph is called {\em simple} if is it composed of distinct vertices.)   
In fact, examining cycle multipliers is the key idea to a better understanding of box expansion and box metrics. We would like to thank Clark Robinson for asking this question.

\begin{defn} \label{defn:avgboxcycmult}
Call $\lambda_k = \min\{\onorm{f'(z)} \colon z \in B_k \}$ the \textit{multiplier} of $B_k$.  
If $B_{0} \to \ldots \to B_{{n-1}} \to B_{n}=B_{0}$ is an $n$-cycle of boxes in $\Gamma$, 
then the {\em cycle multiplier} is $(\lambda_{0} \cdots 
\lambda_{{n-1}})$, and the
{\em average cycle multiplier} is $(\lambda_{0} \cdots 
\lambda_{{n-1}})^{1/n}$.
\end{defn}

Along the way to proving Theorem~\ref{thm:alg}, we establish the following characterizations of box expansion, which are independent of any algorithm used to find an expanded {\boxmet}.  

\begin{prop} 
\label{prop:bestL}
Let $\mathcal L$ be the minimum average cycle multiplier over all simple
cycles in the graph $\Gamma$.  
Then for any $L>0$, there exists a {\boxmet } for which $\Gamma$ is
{\boxexp } by $L$ if and only if  $L \leq \mathcal L$.
\end{prop} 

\begin{cor} \label{cor:minmult} 
Let $\mathcal{M}$ be the minimum cycle multiplier over all simple cycles in the graph $\Gamma$.
Then $\Gamma$ is {\boxexp } (hence $f$ is hyperbolic) if and only if  $\mathcal{M}>1$.
\end{cor}

Thus if some $\Gamma$  is {\boxexp}, then each cycle multiplier for $\Gamma$ is greater than one.  However, some boxes in a cycle could have multiplier less than one, if others are large enough to compensate.  In this case, in the euclidean metric, the map is not expansive along every edge.  The handicaps are designed to spread the expansion out along the cycles, so that in the {\boxmet}, the map is expansive by at least $L$ on every edge in the graph.   This is why we use the term ``handicap''.

According to Corollary~\ref{cor:minmult}, we can show whether $f$ is {\boxexp} simply by computing $\mathcal{M}$. There do exist efficient algorithms for finding such a minimum, see \cite{TCR}.  
However, in order to have explicit information about the hyperbolic structure, we still want to find a viable box expansion constant $L < \mathcal{L}$, and build an expanded {\boxmet}, which we can do  using the {\namealg}.

One weakness of this algorithm is that the box expansion constant, $L>1$, must be inputted in advance.  
Proposition~\ref{prop:bestL} shows the ideal box expansion constant is the smallest average cycle multiplier, $\mathcal{L}$.  
However, in Section~\ref{sec:idealL}, we describe counterexamples which suggest that the only algorithm for explicitly computing $\mathcal{L}$  is exponential (thus too inefficient for our examples).
But then we describe an efficient method for finding a good approximation to the ideal $\mathcal{L}$, in Section~\ref{sec:approxL}. 

Our approach for testing hyperbolicity has some similarities to work of Osipenko (\cite{Osi00,Osi03}).   
For $f$ a diffeomorphism of a compact Riemannian manifold $M \subset \Rn$, he
uses a symbolic image $\Gamma$ of $f$, and develops a general  algorithm for describing the expansion of  $f$ by approximating Lyapunov exponents and the Morse spectrum of the chain recurrent set. His algorithm can also be used to verify hyperbolicity.  Our work here differs in that we are interested in developing and implementing efficient algorithms for families
of polynomial maps of $\CC$ (of degree greater than one).  For this study, we found it more efficient
to get expansion information and a hyperbolicity test by constructing a metric expanded by $f$.  Osipenko's hyperbolicity test has the same computational complexity as finding the smallest average cycle multiplier for all simple cycles in the graph $\Gamma$ (see Section~\ref{sec:idealL}).

In implementation, we control round-off error using {\em interval arithmetic} (IA). This method was recommended
by Warwick Tucker, who used it in his recent computer proof that the
Lorenz differential equation has the conjectured geometry (\cite{War}).
In designing our algorithms, we must keep in mind the workings of IA.  We thus give a brief description of IA in Section~\ref{sec:IA}.


To summarize the organization of the remaining sections: 
in Section~\ref{sec:met1exp} we show box
expansion implies the standard definition of expansion, to prove Theorem~\ref{thm:boxmetexp};
in Section~\ref{sec:IA}, we briefly describe interval arithmetic; 
in Section~\ref{sec:basicalg} we give our basic algorithm, the \namealg, for attempting to 
establish box expansion; 
in Section~\ref{sec:charboxexp} we obtain dynamical information
from either success or 
failure of the {\namealg}, proving Theorem~\ref{thm:alg}, Proposition~\ref{prop:bestL},
and Corollary~\ref{cor:minmult}; 
in Section~\ref{sec:betterL} we compare and contrast ideal versus efficient methods for
determining a good expansion constant $L$; and 
in Section~\ref{sec:examples} we discuss our implementation of the algorithm and 
give examples of output.

\section{Box expansion implies continuous expansion}
\label{sec:met1exp}

In this section, we show box expansion implies 
the standard definition of expansion. Throughout, let $f$ denote a polynomial map of $\CC$ of degree $d>1$, and as in Definition~\ref{defn:boxmodel}, let $\Gamma$ be a {\boxmod } of $f$ on $J$, with vertex set $\mathcal{V}  = \{ B_k \}_{k=1}^N$, composing the {\boxJ } $\boxcov = \cup_{k=1}^N B_k \supset J$.

First, it is more natural for computer calculations, and reduces round-off error, 
to consider vectors in $\Rt$,
rather than $\CC$, and use the $L^{\infty}$ metric of $\Rt$, rather than 
euclidean. Hence, we consider
\[
\sonorm{z}=\max\{ \abs{\text{Re}(z)},\abs{\text{Im}(z)} \}.
\]
Also, let ${\nbd{S}{r}}$ denote the open $r$-neighborhood about the set
$S$ in the metric induced by the above.
This metric is {\em uniformly
equivalent} to the euclidean metric $\eonorm{\cdot}$, since 
$\frac{1}{\sqrt{2}} \eonorm{\cdot} \leq \sonorm{\cdot} \leq \eonorm{\cdot}$.
Thus neighborhoods are slightly different, but the topology generated by them is exactly the same, so they can nearly be used interchangeably.  

To prove Theorem~\ref{thm:boxmetexp}, we use a partition of unity to smooth out a {\boxmet } $\Phi$.  
First, we need a lemma from \cite{SLHone}.   

\begin{lem}[\cite{SLHone}] \label{lem:imageboxsize}
There exists an $\eta>0$ so that if $B_k, B_j \in \mathcal{V}$, with $z \in \nbd{B_k}{\eta}$ and $f(z) \in \nbd{B_j}{\eta}$, then there is an edge from $B_k$ to $B_j$ in $\Gamma$.
\end{lem}

To prove this lemma, we used the assumption that $f$ was a polynomial of degree $d > 1$, and the fact that by Definition~\ref{defn:boxmodel}, there is a $\delta>0$ such that 
there is an edge from $B_k$ to $B_j$ 
if a $\delta$-neighborhood of $f(B_k)$ intersects $B_j$.  Now, we can obtain:

\begin{lem} \label{lem:partnbd}
Suppose $\Gamma$ is \boxexp.
Then there exists
a $\tau >0$ such that if $B_k,B_j \in \mathcal{V}$, and
$z\in \boxcov \cap \nbd{B_k}{\tau}$ 
with $f(z) \in \nbd{B_j}{\tau}$, then for any $\mathbf{v}\in T_z \CC$,
$\jonorm{D_zf(\mathbf{v})} \geq L \konorm{\mathbf{v}}$.
\end{lem}

\begin{proof}
Among other requirements below, let $\tau>0$ be less than $\eta$ from 
Lemma~\ref{lem:imageboxsize}.   Then for $z$ satisfying the hypotheses,
there is an edge from $B_k$ to $B_j$ in $\Gamma$.

Note since we are working in one dimension, $D_x f=f'(x)$, hence box expansion
yields that for $x\in B_k$,
$\varphi_j \abs{f'(x)} \abs{\mathbf{v}} \geq L \varphi_k \abs{\mathbf{v}}$, thus 
simply
$\varphi_j \abs{f'(x)} \geq L \varphi_k$.

Since $\mathcal B$ is compact, $\mathcal{V}$ is finite, and $f'(x)$ is continuous,
there is an $\alpha \geq 0$ so that:
\begin{enumerate}
\item  $\alpha = \min \{ \varphi_j \abs{f'(x)} - L\varphi_k \colon x\in B_k, (k,j) \in 
\mathcal{E} \}$,
\item  if $\tau <\eta$ is sufficiently small, then for any $j$,
$\abs{x-z} < \tau$ implies that
 $\varphi_j \abs{f'(x)-f'(z)} < \alpha$.
\end{enumerate}

Now $z$ is not necessarily in $B_k$, but $z\in \boxcov$, so suppose
$z\in B_m$ and $x\in B_m \cap B_k$ such that $\abs{x-z} < \tau$.  Then
$\varphi_j \abs{f'(x)-f'(z)} < \alpha$; further, there is an edge
$(k,j)\in\mathcal{E}$, hence $x$ satisfies $\varphi_j \abs{f'(x)} - L\varphi_k \geq
\alpha$. Combining these gives $\varphi_j \abs{f'(z)} \geq L \varphi_k$. Thus
$\jonorm{D_zf(\mathbf{v})} \geq L \konorm{\mathbf{v}}$.
\end{proof}

Now we use the $\tau$-overlap to convert a {\boxmet } into a riemannian metric.

\begin{defn} \label{defn:ponorm}
Suppose $\Gamma$ is {\boxexp } for a {\boxmet } $\Phi$.
Let $\tau >0$ be as given by Lemma~\ref{lem:partnbd}. 
Define a partition of  unity on $\boxcov(\Gamma)$ by
choosing smooth
 functions $\rho_k:\CC \to [0,1]$, for each box 
$B_k \in \mathcal{V}$, 
such that $\text{supp}(\rho_k) \subset \nbd{B_k}{\tau}$
and $\sum_{k} \rho_k(x) = 1$, for any $x \in \mathcal B$.
Define the smooth function 
$\rho = \rho(\Phi)  \colon \boxcov \to [0,1]$ by $\rho(x) = \sum_{k} \rho_k (x) \varphi_k$. Then
$\rho$ induces a riemannian metric on $T\boxcov$,
with a smoothly varying norm, $\ponorm{\cdot}$, where if $x\in \boxcov$
and $\mathbf{v} \in T_x \boxcov$ then 
\[
\ponorm{\mathbf{v}} = 
\onorm{\mathbf{v}} \rho(x) =
\onorm{\mathbf{v}} \sum_{k} \rho_k (x) \varphi_k =
\sum_{k}\rho_k (x) \onorm{\mathbf{v}}_k.
\]
\end{defn}

Note $\ponorm{\cdot}$ and $\Phionorm{\cdot}$ are very close.  
They only differ in the small
$\tau$-neighborhoods of the box boundaries, 
where the $\rho$ metric smooths out the $\Phi$ metric.  

It is also straightforward to show that both the 
$\rho$ metric and the $\Phi$ metric are uniformly equivalent to euclidean,
with
\[
\(\min_k \{ \varphi_k \} \)  \onorm{\cdot} \leq \ponorm{\cdot}, \Phionorm{\cdot} \leq   \(\max_k \{ \varphi_k \}\) \onorm{\cdot}.
\]

We establish Theorem~\ref{thm:boxmetexp} by showing that if $\Gamma$ is {\boxexp} for a {\boxmet } $\Phi$, then $f$ is expansive in $\boxcov$ for the metric $\rho(\Phi)$.

\begin{proof} [Proof of Theorem~\ref{thm:boxmetexp}]
Let $\ponorm{\cdot}, \tau$ be as in Definition~\ref{defn:ponorm}.  

Thus $\tau$ is small enough that if $x \in \text{supp} 
(\rho_k)$, then for any $j$ such that $f(x) \in \text{supp}(\rho_j)$,
we have
$
\varphi_j\abs{D_xf(\mathbf{v})} \geq L \varphi_k \abs{\mathbf{v}},
$
for any $\mathbf{v} \in T_x\CC$.
Then if we set
\[
\varphi_x = \max\{ \varphi_k \colon x \in \text{supp}(\rho_k) \}, \text{ \ and \ }
\varphi_{f,x} = \min \{ \varphi_j \colon f(x) \in \text{supp}(\rho_j) \},
\]
we know 
$
\varphi_{f,x} \abs{D_xf(\mathbf{v})} \geq L \varphi_x \abs{\mathbf{v}},
$
for any $\mathbf{v} \in T_x\CC$.

Now if we use that $\sum_{k} \rho_k(x) = 
\sum_{j} \rho_j(f(x)) = 1$, we get the 
result:
\begin{eqnarray*}
\ponorm{D_xf(\mathbf{v})} &=& \abs{D_xf(\mathbf{v})} \sum_{j} \rho_j(f(x)) \varphi_j
\geq \abs{D_xf(\mathbf{v})} \varphi_{f,x} \\
& \geq & L  \varphi_x \abs{\mathbf{v}} 
\geq L \abs{\mathbf{v}} \sum_{k} \rho_k(x) \varphi_k
= L \ponorm{\mathbf{v}}.
\end{eqnarray*}
\end{proof}

Thus if any box model $\Gamma$ of $f$ on $J$ 
 is {\boxexp } by some $L>1$, then $f$
 is expanding on $J$ in the riemannian metric 
$\rho$, and hence $f$ is hyperbolic.

\section{Interval Arithmetic}
\label{sec:IA}

Interval arithmetic (IA) provides a natural and efficient method for manipulating boxes,
and also for maintaining rigor in computations.  The basic objects of IA
 are closed intervals, $[a] = [\underline{a}, \bar{a}] \in \II\KK$, with end points in some fixed
 field, $\KK$.   An arithmetic operation on two intervals produces a resulting interval
which contains the real answer.  For example,
\begin{align*}
[a] + [b] & := 
\left[ \underline{a}+\underline{b}, \bar{a}+\bar{b} \right] \\
[a] - [b] & := 
\left[ \underline{a}-\bar{b}, \bar{a}-\underline{b} \right] 
\end{align*}
Multiplication and division can also be defined in IA.

Computer arithmetic is performed not with real numbers, but rather in
the finite space $\FF$ of numbers representable by binary floating point
numbers of a certain finite length.  For example, since the number 
$0.1$ is not a dyadic rational, it has an infinite binary expansion,
so is not in $\FF$.   

Since an arithmetical operation on two numbers in $\FF$ may not have a 
result in $\FF$,  in order to implement rigorous IA we must round
outward the result of any interval
arithmetic operation, \textit{e.g.} for $[a],[b] \in \II\FF$,
$$
[a] + [b]  := 
\left[ \left\downarrow \underline{a}+\underline{b} \right\downarrow, 
         \left\uparrow  \bar{a}+\bar{b} \right\uparrow \right],
$$
where $\left\downarrow x \right\downarrow$ denotes the largest number in $\FF$ 
that is 
strictly less than $x$ (\textit{i.e.}, $x$ rounded down), and 
 $\left\uparrow x 
\right\uparrow$ 
denotes the smallest number in $\FF$ that is strictly greater than $x$ 
(\textit{i.e.}, $x$ 
rounded up).  This is called IA with \textit{directed rounding}.

For any $x\in \RR$, let Hull$(x)$ be the smallest interval in $\FF$ which contains $x$.
That is, if $x \in \FF$, then Hull$(x)$ denotes $[x, x]$.  If $x \in \RR \setminus \FF$, then 
 Hull$(x)$ denotes $\left[ \left\downarrow x \right\downarrow, 
\left\uparrow x \right\uparrow \right]$.

In higher dimensions, IA operations can be carried out component-wise, on 
{\em interval vectors}. 
Note a box in $\CC = \Rt$ is simply
an interval vector.

In designing our algorithms, every arithmetical calculation is carried out with IA, and thus we must think carefully
about how to use IA in each situation.  
Our extensive use of boxes is designed to 
make IA calculations natural.
However, IA tends to create problems with
propagating increasingly large error bounds, if not handled carefully.
For example, iterating a polynomial map on an interval which is close to the Julia set, $J$, can produce a tremendously large interval after only a few iterates (due to the expanding behavior of the map 
near $J$).  
Thus, in the remainder of the paper, after describing each algorithm we note how IA is being used.
\begin{footnote}
{All of our IA computations use the PROFIL/BIAS package, available at \cite{PBIA}.}  
\end{footnote}
The interested reader can find an abundance of materials on IA, for example \cite{GenIA, MooreIA1, MooreIA2}.

\section{The {\namealg }}
\label{sec:basicalg}

In this section we describe our basic algorithm, the {\namealg}, for attempting to build an expanded {\boxmet } for a {\boxmod } $\Gamma$.

The problem of finding a set of handicaps $ \{ 
\varphi_k \}$ defining a {\boxmet } $\Phi$ for which $\Gamma$ is {\boxexp } by a given $L>1$ is strictly a graph theoretic problem.  
We want to find handicaps so that for every edge $(k,j) \in 
\mathcal{E}$, for every point $z \in B_k$, and for every vector $\mathbf{v} \in T_z 
B_k$, 
we have
\[
\varphi_j\sonorm{D_zf(\mathbf{v})}
= \varphi_j \sonorm{D_zf} \sonorm{\mathbf{v}}
\geq L \varphi_k  \sonorm{\mathbf{v}},
\]
or equivalently, if $\lambda_k = \min_{z \in B_k} \sonorm{f'(z)}$  is the multiplier of $B_k$,
we want
\begin{equation} \label{eqn:metricfactor}
\varphi_j \geq 
L \varphi_k / \lambda_k.
\end{equation}

For polynomial maps of $\CC$, a lower bound for the
multiplier $\lambda_k$ in a box is easily calculated with interval arithmetic.
After the $\lambda_k$ are calculated, we can forget the map 
$f$ and simply think of $\Gamma$ as a strongly connected directed graph 
endowed with edge
weights, $\xi_{k,j} = L/ \lambda_k$, if $(k,j)\in\mathcal{E}$.  In the case of maps in one dimension, the edge weights are all equal along all edges emanating from a single vertex, \textit{i.e.}, $\xi_{k,1} = \xi_{k,2} = \cdots \xi_{k,n}$.  However we describe the algorithm for the more general situation, where $\xi_{k,j} \neq \xi_{k,m},$ since it is not more difficult, and it is useful for higher dimensional applications (for example, in \cite{SLHthree}). 

\begin{defn} \label{defn:handicap}
Let $\Gamma$ be a directed graph with positive edge 
weights $\Xi=$ $ \{\xi_{k,j} \colon$ $ (k,j)\in \mathcal{E} \}$.  
Let $\Phi = $ $ \{\varphi_k \colon $ $v_k\in\mathcal{V} \}$ be a set of positive 
vertex weights, called \textit{handicaps}.  We call $\Phi$ {\em 
consistent handicaps} (for $\Gamma, \Xi$) if 
$\varphi_j \geq \xi_{k,j} \ \varphi_k$, for every edge $(k,j)\in 
\mathcal{E}$.  

If equality holds on every edge, then we call $\Phi$ {\em 
strict handicaps}.
\end{defn}

Then $\Gamma$ is {\boxexp } if there 
exists a set of consistent handicaps $\Phi$ for $\Gamma$, given the 
edge weights suggested by~(\ref{eqn:metricfactor}).   

To attempt to find consistent handicaps given $\Gamma$ and the edge weights $\Xi$, 
we break
the problem into a finite induction. Step $0$ consists of finding an
initial set of strict handicaps for some spanning tree, $\Gamma_0$, of
$\Gamma$.  At steps $n>0$, we choose a graph $\Gamma_n$ such that
$\Gamma_0 \subset \Gamma_n \subset \Gamma$, then seek a consistent
set of handicaps for $\Gamma_n$.

\begin{defn} \label{defn:arbor}
A directed graph $\Gamma_0$ is an {\em arborescence} if there is 
a root vertex $v_0$ so that for any other vertex $v$, 
there is a unique simple path from $v_0$ to $v$.  Such a graph is a 
tree, and must have exactly one incoming edge for each vertex $v \neq 
v_0$.  

If $\Gamma$ is strongly connected, then for each vertex 
$v_0$ in $\Gamma$, there is a minimum spanning tree $\Gamma_0$ with 
root vertex $v_0$ which is an arborescence 
(simply perform a depth first or breadth first search from $v_0$).  We 
call such a $\Gamma_0$ a {\em spanning arborescence}.
\end{defn}

\begin{defn}
Let $\Gamma$ be a finite, strongly connected directed graph. Let $\Gamma_0 \subset 
\Gamma_1 \subset \cdots \subset \Gamma_N = \Gamma$
be a nested sequence of subgraphs of $\Gamma$ such that
$\Gamma_0$ is a spanning arborescence of $\Gamma$, and for $n\geq 1$,
$\mathcal{V}(\Gamma_n) = \mathcal{V}(\Gamma)$, and $\mathcal{E}(\Gamma_n)$ is
formed by adding one edge of $\mathcal{E}(\Gamma) \setminus \mathcal{E}(\Gamma_{n-1})$
to $\mathcal{E}(\Gamma_{n-1})$.  We call such a
sequence an {\em edge exhaustion of $\Gamma$}.
\end{defn}

Note by the above definitions that an edge exhaustion exists for any finite, strongly connected  directed graph $\Gamma$.
Note also that since $\Gamma_0$ is a spanning arborescence, each $\Gamma_n$ is edge connected.

\begin{alg}[Recursively hedging handicaps via an edge exhaustion]
\label{alg:hyponemetric}
Let $\Gamma_0 \subset \Gamma_1 \subset \cdots \subset \Gamma_N = \Gamma$ 
be an edge exhaustion of $\Gamma$.

\textbf{Base Case:}
Inductively construct a set of {\em strict} handicaps, $\Phi^0 = \{ \varphi^0_k \}$ for
$\Gamma_0$, by choosing any $\varphi^0_0$ for the root vertex $v_0$, then pushing this
value across $\Gamma_0$ by multiplication with successive edge weights.  
That is, if $v_k \in \mathcal{V}(\Gamma_0)$, and $v \neq v_0$, then since $\Gamma_0$ is an arborescence, the predecessor, $v_{\pi(k)}$, is uniquely defined such that the edge $(\pi(k),k)$ is in $\Gamma_0$.
Then set $\varphi^0_k = \xi_{ k, \pi(k)} \ \varphi^0_{\pi(k)}$.

For our box expansion application, we pick a random vertex as the
root, and start with $\varphi^0_0=1$, \textit{i.e.}, take the euclidean metric 
on that box.

\textbf{Inductive Step:}
Suppose we have a set of consistent handicaps $\Phi^{n-1} = \{ \varphi^{n-1}_k \}$ for $\Gamma_{n-1}$.  
We attempt to define a set of consistent handicaps $\Phi^n = \{ \varphi^n_k \}$ for $\Gamma_{n}$, by
adjusting the set $\Phi^{n-1}$, using a process we call {\em 
recursively hedging the handicaps}.  First, start with the set of temporary handicaps $\Phi^n := \Phi^{n-1}$ on $\Gamma_n$.
Let $(v_a,v_b)$ denote the edge in $\Gamma_{n}$ 
that is not in $\Gamma_{n-1}$ (this edge is unique by definition of edge exhaustion). Then we adjust the $b$th handicap, by doing
\begin{tabbing}
\hspace*{.2in} \= \hspace{.2in} \= \hspace{.2in}
\kill
\>if $\varphi^n_b < \xi_{a,b} \ \varphi^n_a $\\
\>\>set $\varphi^n_b = \xi_{a,b} \ \varphi^n_a $
\end{tabbing}

If the second line above is performed, \textit{i.e.} $\varphi^n_b$ is strictly increased, then we say we have {\em hedged} $\varphi^n_b$ {\em along the edge} $(v_a, v_b)$.  But now if $\varphi^n_b$ is hedged, then  
in order to keep consistency of the handicaps, for each vertex $v_c$ in $\Gamma_n$ which is
adjacent to $v_b$, we may need to hedge $\varphi^n_c$ along the edge $(v_b, v_c)$:
\begin{tabbing}
\hspace*{.2in} \= \hspace{.2in} \= \hspace{.2in} \= \hspace{.2in} \= \hspace{.2in}
\kill
\>\>set $\varphi^n_b = \xi_{a,b} \ \varphi^n_a $\\
\>\>for each $v_c$ in $\Gamma_n$ and in the adjacency list of $v_b$ do\\
\>\>\>if $\varphi^n_c < \xi_{b,c} \ \varphi^n_b $\\
\>\>\>\>set $\varphi^n_c = \xi_{b,c} \ \varphi^n_b $
\end{tabbing}

The hedging process becomes recursive here, for if any $\varphi^n_c$ is increased,
then we may need to hedge handicaps along each edge in $\Gamma_n$ emanating from
$v_c$, etc.  The worst case scenario is that every vertex reachable from $v_b$ along every
path in $\Gamma_n$ would need to be hedged.

However if at any step in the hedging, we find a vertex $v_x$ whose handicap need not be
increased, \textit{i.e.}, $\varphi^n_x \geq \xi_{\pi(x),x} \ \varphi^n_{\pi(x)}$ already, then we may stop, since the handicaps of vertices reachable from $v_x$ will not need to be increased. This greatly saves on computational time.
Thus, our process so far is:
\begin{tabbing}
\hspace*{.2in} \= \hspace{.2in} \= \hspace{.2in} \= \hspace{.2in} \= \hspace{.2in}
\kill
\>\textbf{Hedge$(a,b)=$}\\
\>if $\varphi^n_b < \xi_{a,b} \ \varphi^n_a $\\
\>\>set $\varphi^n_b = \xi_{a,b} \ \varphi^n_a $\\
\>\>for each $v_c$ in $\Gamma_n$ and in the adjacency list of $v_b$ do\\
\>\>\>\textbf{Hedge($b,c$)}
\end{tabbing}

Will this recursive procedure always terminate in a consistent set of
handicaps $\Phi^n$ for $\Gamma_n$?  Of course not, since it is possible
that there does not exist any set of consistent handicaps for $\Gamma$
with edge weights $\Xi$.
There is one possible obstruction: suppose that for the new edge in 
$\Gamma_n$, edge $(v_a,v_b)$, we travel away from $v_b$, hedging the handicaps 
along every edge in some path, $v_a, v_b, v_c, \ldots, v_x$,
and then discover $v_a$ in the adjacency list of $v_x$.  This means every
handicap on a path in $\Gamma_n$ from $v_b$ to $v_x$ has been
increased, since we stopped the search if it was not.  The handicap
at $v_a$ cannot be increased, since then we would have to hedge $\varphi^n_b$
along $(v_a,v_b)$ again, and an infinite loop of hedging along this path
through $v_x$ would occur. 

Thus the cycle $\{v_a, v_b, v_c, \ldots, v_x, v_a \}$ 
could be an obstruction to consistent handicaps,
for upon seeing $v_a$ from $v_x$, we can only check if $\varphi^n_a \geq 
\xi_{x,a} \ \varphi^n_x$ already.  This leads us to the following algorithm for
the $n^{th}$ inductive step: 

\begin{tabbing}
\hspace*{.2in} \= \hspace{.2in} \= \hspace{.2in} \= \hspace{.2in} \= \hspace{.2in}
\kill
\>do \textbf{Hedge$(a,b,a)$}, where\\
\>\textbf{Hedge$(u,v,a)=$}\\
\>if $\varphi^n_v \geq \xi_{u,v} \ \varphi^n_u $, then return 1\\
\>else if $(v=a)$, then return 0\\
\>else \\
\>\>set $\varphi^n_v = \xi_{u,v} \ \varphi^n_u $\\
\>\>for each $w$ in $\Gamma_n$ and in adjacency list of $v$ do\\
\>\>\>\textbf{Hedge$(v,w,a)$}
\end{tabbing}

Thus the $n^{th}$ step is terminated when either a cyclic
obstruction is found (at the ``return 0'' line above), or when all
necessary hedgings have successfully been performed
(if  \textbf{Hedge$(a,b,a)$} returns 1).  If no cyclic obstructions
are found, then $\Phi^n$ is a consistent set of
handicaps.
\end{alg}

Thus Algorithm~\ref{alg:hyponemetric} either terminates in a cyclic obstruction at
some $\Gamma_n$, or produces consistent handicaps for all of
$\Gamma=\Gamma_N$. This dichotomy leads to Theorem~\ref{thm:alg}, proved
in Section~\ref{sec:charboxexp}.


The {\namealg } consists simply of
applying Algorithm~\ref{alg:hyponemetric} to the case of showing box expansion for a {\boxmod } $\Gamma$ of $f$ on $J$.
We describe this process more explicitly with the following pseduo-code.
Comments are parenthetical and to the right.
The routine consists of a main function, \textbf{BuildMetric}, and its two helper sub-routines, \textbf{SpanTree} and \textbf{Hedge}.   
\textbf{SpanTree} creates an arborescence using a breadth first search and a queue (first-in $=$ first-out) to traverse the graph.  \textbf{Hedge} uses a depth first style to check the rest of the edges.  See \cite{TCR} for background on graph searches, and note the remark following the pseudo-code on why these styles were chosen.

\begin{alg}[The \namealg]
 \label{alg:fexpands}
\begin{tabbing}
\hspace*{.2in} \= \hspace{.2in} \= \hspace{.2in} \= \hspace{.2in}
\= \hspace{.2in} \= \hspace{.2in} \= \hspace{.2in} \= \hspace{.2in} \= \hspace{.2in}
\= \hspace{.2in} \= \hspace{.2in} \= \hspace{.2in} \= \hspace{.2in} \= \hspace{.2in}
\kill
\>\>\\
\> \textbf{BuildMetric($\Gamma, L$):}\\
\>for every vertex $u$ in $\Gamma$\\
\>\>set color$[u]$ to white\\
\>\>for each vertex $v$ in adjacency list of $u$ do\\
\>\>\>set check$[u][v]$ = 0\\
\>set $\varphi_0 = 1$\\
\>do \textbf{SpanTree}(vertex $0$ of $\Gamma$)\\
\>for every edge $(a, b)$ in $\Gamma$\\
\>\>if (check$[a][b] $= 0) then 
	\>\>\>\>\>\>\>{(\textit{edge $(a,b)$ not in previous step of exhaustion})}\\
\>\>\>set check$[a][b]$ = 1
	\>\>\>\>\>\>\>{(\textit{put $(a,b)$ into current step of exhaustion})}\\
\>\>\>if (\textbf{Hedge}$(a,b,a,L) = 0)$ then return $0$\\
\>return $1$  
 \>\>\>\>\>\>\>\>\>(\textit{if get this far, then successful})\\

\>\>\\
\>\textbf{SpanTree($u$):}\\
\>put vertex $u$ on the queue, $Q$ (at the end)\\
\>set color[$u$] to gray\\
\>while $Q \neq \emptyset$\\
\>\>let $v$ be the head (first element) of the queue\\
\>\>for each vertex $w$ in adjacency list of $v$ do\\
\>\>\>if color[$w$] is white\\
\>\>\>\>put vertex $w$ on the queue, $Q$ (at the end)\\
\>\>\>\>set color[$w$] to gray\\
\>\>\>\>set check[$v$][$w$] = 1
\>\>\>\>\>\>(\textit{record edge $(v,w)$ is in spanning tree})\\
\>\>\>\>set $\varphi_w = L \varphi_v / \lambda_v$\\
\>\>remove first element from queue 
\>\>\>\>\>\>\>\>\>(\textit{remove $v$})\\
\>\>set color$[v]$  to black 
\>\>\\
\>\>\\
\>\textbf{Hedge($u$, $v$, $a$, $L$):}\\
\>if ($\varphi_v \geq L \varphi_u/ \lambda_u$) then return $1$ 
\>\>\>\>\>\>\>\>\>\>(\textit{i.e., edge $(u,v)$ is already ok})\\
\>else if (v = a) then return $0$
\>\>\>\>\>\>\>\>\>\>(\textit{we cannot increase $\varphi_a$, so fails})\\
\>else\\
\>\>set $\varphi_v = L \varphi_u/ \lambda_u$ 
 \>\>\>\>\>\>\>\>\> (\textit{increase $\varphi_v$})\\
\>\>for each $w$ in adjacency list of $v$ do\\
\>\>\>if (check[v][w] = 1) then
 \>\>\>\>\>\>\>(\textit{$(v,w)$ is in current step of exhaustion}) \\
\>\>\>\>if (\textbf{Hedge}$(v,w,a,L) = 0)$ then return $0$\\
\>return $1$ 
\>\> \>\>\>\>\>\>\>\>(\textit{if get this far, then successful})
\end{tabbing}
\end{alg}

\begin{rem}
In \textbf{SpanTree} we maximize efficiency by using a breadth first style search to traverse the graph, instead of depth first.
Running depth first search on a typical {\boxmod } $\Gamma$ tends to produces a spanning arborescence with a very long path, whereas breadth first search constructs an arborescence with paths of minimum length.  This is because for a polynomial map $f$, the typical $\Gamma$ has a large number of vertices and a small bound on the out-degree of vertices (related to the bound on the derivative of the map in $\boxcov$).  In addition, $\Gamma$ is strongly connected.  

Long paths creates two problems. First, creating such a tree can causes a memory overflow, due to  many nested recursive function calls.  Additionally, even when there is no crash, a large path tends to create an initial metric farther from Euclidean. See Example~\ref{exmp:smallL} in
Section~\ref{sec:betterL} for an exploration of the latter phenomena.
\end{rem}

Recall in our implementation, we control rounding using interval arithmetic (Section~\ref{sec:IA}).   In the above, we round up to ensure the inequality is satisfied, \textit{i.e.}, check
$
\varphi_v \geq \sup( \text{Hull}(L) *\text{Hull}(\varphi_u)/ \text{Hull}(\lambda_u)).
$

\section{Characterization of box expansion}
\label{sec:charboxexp}

Recall from Definition~\ref{defn:avgboxcycmult} that the {multiplier} of a box $B_k$ is $\lambda_k = \min_{z\in B_k} \abs{f'(z)}$,
and for an $n$-cycle of boxes $(B_{0}\to \ldots \to B_{{n-1}} \to B_{n}=B_{0})$ in $\Gamma$, we call $(\lambda_{0} \cdots 
\lambda_{{n-1}})$ the {cycle multiplier} and $(\lambda_{0} \cdots 
\lambda_{{n-1}})^{1/n}$ the {average cycle multiplier}. Also, a cycle in a graph is called {\em simple} if is it composed of distinct
vertices.  

We can now specify
the implications of success or failure in the \namealg,  
proving Theorem~\ref{thm:alg}, Proposition~\ref{prop:bestL}, and Corollary~\ref{cor:minmult}.

\begin{prop} \label{prop:hyponefailed}
For $L>0$, the \namealg, \textit{i.e.}, Algorithm~\ref{alg:fexpands}, or
Algorithm~\ref{alg:hyponemetric} with edge weights $\Xi = \{ \xi_{k,j} = L/\lambda_k \colon
(k,j) \in \mathcal{E}\}$, 
either constructs a set of consistent handicaps showing $\Gamma$ is {\boxexp } by $L$, 
or finds an $n$-cycle of boxes with cycle multiplier  less than $L^n$.
\end{prop}

\begin{proof}
We observed in the description of Algorithm~\ref{alg:hyponemetric} that 
the only obstruction to building consistent handicaps
is if in adding an edge $(u, v)$, we find a cycle of boxes
$(u=B_{0} \to v=B_{1} \to \ldots \to B_{{n-1}} \to B_{n} = B_{0})$, such that,
holding $\varphi_{0}$ fixed, in order to keep consistency the metric
handicaps must be increased along every edge in the cycle. That is,
$\varphi_{{k+1}} = L \varphi_{k} / \lambda_{k}$ for $0 \leq k \leq n-2$,
and we have the failure $\varphi_{0} < L \varphi_{{n-1}} /
\lambda_{{n-1}}$.  But then
\[
\varphi_{0} < 
\frac{L \varphi_{{n-1}}}{\lambda_{{n-1}}} =
\frac{L^2 \varphi_{{n-2}}}{\lambda_{{n-1}} \lambda_{{n-2}}} =
\cdots =
\frac{L^n \varphi_{0}}{\lambda_{{n-1}} \lambda_{{n-2}} \cdots
\lambda_{{0}}}.
\]
Hence, 
$\lambda_{{n-1}} \lambda_{{n-2}}\cdots\lambda_{{0}} <
L^n$.
\end{proof}

\begin{lem} \label{lem:hyponeworked}
Let $(B_{0} \to \ldots \to B_{{n-1}} \to B_{n} = B_{0})$ be an $n$-cycle of boxes in
$\Gamma$, such that consistent handicaps $\{\varphi_{0}, \ldots,
\varphi_{{n-1}}, \varphi_{n} = \varphi_{0}\}$ can be chosen to show $f$ box expands by $L>1$
along the cycle.  
Then its cycle multiplier is at least $L^n$.
\end{lem}

\begin{proof}
By hypothesis, we know
$
\varphi_{{k+1}} \lambda_k \geq L \varphi_{k}, \
\ k \in \{0,\ldots,n-1\}.
$  
Thus,
\[
 \lambda_0 \cdots \lambda_{n-1} \geq 
 L \frac{\varphi_{0}}{\varphi_{1}} L \frac{\varphi_{1}}{\varphi_{2}} \cdots 
  L \frac{\varphi_{{n-1}}}{\varphi_{n} = \varphi_{0}}.
\]
Cross cancellation and simplifying leaves only
$
\lambda_0 \cdots \lambda_{n-1} \geq
L^n.
$
\end{proof}

\begin{lem} \label{lem:failsiffcycle}
For $L>0$, there is a cycle of boxes in $\Gamma$ with average multiplier less than $L$ if and only if there is no {\boxmet } which $\Gamma$ expands by $L$.
\end{lem}

\begin{proof}
First, Lemma~\ref{lem:hyponeworked} gives the forward implication
immediately;
for, if there is a {\boxmet } which $\Gamma$ expands by $L$,
then every cycle
must have average multiplier at least $L$.
Conversely, suppose there is no {\boxmet } which
$\Gamma$ expands by $L$.   
Then by Proposition~\ref{prop:hyponefailed}, there exists a 
cycle with average multiplier less than $L$.
\end{proof}

\begin{proof}[Proof of Theorem~\ref{thm:alg}]
Proposition~\ref{prop:hyponefailed} shows 
that given $L>0$, the {\namealg } either constructs a {\boxmet } for which $\Gamma$ is expansive by $L$, or finds a cycle with average multiplier less than $L$.
Lemma~\ref{lem:failsiffcycle} shows the latter is equivalent to saying there is no {\boxmet } which $\Gamma$ expands by $L$.
\end{proof}

\begin{proof}[Proof of Proposition~\ref{prop:bestL}]
This follows directly from Lemma~\ref{lem:failsiffcycle}.
Indeed, let $L>0$ be given.   
First, assume $L \leq \mathcal{L}$.  Then all cycles have average multiplier
at least $L$, so by Lemma~\ref{lem:failsiffcycle} there exists a {\boxmet } 
which $\Gamma$ expands by $L$.
Next, assume $L > \mathcal{L}$. Then the cycle of boxes with average
multiplier $\mathcal{L}$ has average multiplier less than $L$, hence
by Lemma~\ref{lem:failsiffcycle} there is no {\boxmet } which
$\Gamma$ expands by $L$.
\end{proof}

\begin{proof}[Proof of Corollary~\ref{cor:minmult}]
This follows easily from Proposition~\ref{prop:bestL}, using the fact that $a>1$ if and only if $a^{1/m} >1$, for any positive $a$ and $m$.  
Let $\mathcal{L}$ denote the smallest average cycle multiplier over all simple cycles in the graph.
By Proposition~\ref{prop:bestL},   
$\Gamma$ is box expansive if and only if $\mathcal{L} \geq 1$. 
Let $m$ be the length of the cycle with multiplier $\mathcal{M}$.
Let $n$ be the length of the cycle with average multiplier $\mathcal{L}$.
Then $\mathcal{M}^{1/m} \geq \mathcal{L}$ and $\mathcal{M} \leq \mathcal{L}^n$.
Suppose $\Gamma$ is {\boxexp}.   Then $\mathcal{M}^{1/m} \geq \mathcal{L} > 1$, hence $\mathcal{M} > 1$.
Conversely, suppose $\mathcal{M} > 1$.  Then $\mathcal{L}^n \geq \mathcal{M} > 1$,
hence $\mathcal{L} >1$. Thus $\Gamma$ is {\boxexp}.
\end{proof}

\section{Finding a good box expansion constant $L$}
\label{sec:betterL}

%
In theory and in practice, a box expansion constant more appropriate for a particular collection of boxes yields a ``better'' metric.  In particular, in running the program Hypatia it is easy to find behavior
like the following. Suppose for some map, we successfully build a metric
with box expansion constant $L=1.2$, but the resulting  handicaps range between $4.3
\times 10^{-18}$ and $0.6$.  Depending on the map $f$, trying $L =1.5$ could improve the handicaps to be in the range $0.068$ to $0.55$.  This yields a metric more computationally tractable, and closer to euclidean.
Why do we see this behavior?  Consider the example:

\begin{exmp}  \label{exmp:smallL}
Suppose there is a path of boxes $B_{0} \to \ldots \to B_{n}$ such
that (for simplicity) all multipliers are
approximately the same value $\lambda$.  Then suppose we try to
 show box expansion by some $L < \lambda$.  We would 
define the handicaps in the boxes of this cycle so that $\varphi_{{k+1}} \geq
L \varphi_{k} / \lambda$, for $0 \leq k \leq n-1$.  Thus we would get
\[
\varphi_{1} = \frac{\varphi_{0} L}{\lambda},
\varphi_{2} = \frac{\varphi_{1} L}{\lambda}
        = \frac{\varphi_{0} L^2}{\lambda^2},
\ldots,
\varphi_{{n}} = \frac{\varphi_{0} L^n}{\lambda^n}.
\]
But since $L < \lambda$, large $n$ leads to $L^n \ll
\lambda^n$ and $\varphi_{n} \ll \varphi_{0}$.  Thus, if we use an $L$ which is ``too small'', then the handicaps
plummet.
\end{exmp}

In fact in running Hypatia, we have found that 
it is easy for a bad choice of~$L$ to lead to handicaps which are so
large or small that the machine cannot distinguish them from $0$ or
$\infty$.  Checks must be put in place in the algorithms to flag
such occurrences.  
This can be guarded against somewhat by using breadth first 
search style algorithms as much as possible, rather than 
depth first search.  In particular, in \textbf{SpanTree}, if a depth first search is used,
 then the spanning arborescences usually contain very long paths.  
On longer paths, as Example~\ref{exmp:smallL} illustrates,
 the handicaps are more likely to get unmanageable.

\subsection{An optimal, yet impractical solution}
\label{sec:idealL}

Proposition~\ref{prop:bestL} implies that the ideal box expansion constant, $\mathcal{L}$, is equal to the smallest average cycle multiplier over all simple cycles.  This can \textit{theoretically} be computed, since the graph is finite.  However, in this section we shall see the computation seems unwieldy.

For a first approach to finding $\mathcal{L}$, recall that there do exist efficient algorithms for finding the cycle with the smallest cycle multiplier (\cite{TCR}). 
Unfortunately, the following example shows that the smallest average cycle multiplier cannot be simply derived from examining the cycle with the smallest cycle multiplier.

\begin{exmp} \label{exmp:avgbad}
Consider the graph with three vertices $\{0,1,2\}$, and with the four
edges $\{ (0,1), (1,2), (0,2), (2,0) \}$, shown in
Figure~\ref{fig:graphone}.  It has two cycles: $0 \to 1 \to 2 \to 0$ and
$0 \to 2 \to 0$, with multipliers: $\lambda_0 = 1/2,
\lambda_1 = 3/2, \lambda_2 = 8$.  
\begin{figure}
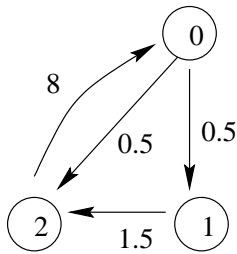

\begin{center}
\drawfiggraphzero
\caption{Graph illustrating Examples~\ref{exmp:avgbad} and~\ref{exmp:avgpathbad}}
\label{fig:graphone}
\end{center}
\end{figure}
The cycle multipliers  are: $(\lambda_0 \lambda_2) = {4}$  and
$(\lambda_0 \lambda_1 \lambda_2) = 6$.  Thus the smallest
cycle multiplier is along the cycle $0 \to 2 \to 0$.  The average along this cycle is $4^{1/2}=2$.
However, the average cycle multiplier along the other cycle is  $(\lambda_0 \lambda_1 \lambda_2)^{1/3} = 6^{1/3}=1.82$.  
\end{exmp}


Thus, we need an algorithm for computing the smallest average cycle multiplier from scratch.  For our second approach to this computation, we might posit that we could adapt the algorithm used to compute the smallest cycle multiplier, to instead compute the smallest average cycle multiplier.  This algorithm involves considering ``best paths'' during a breadth first search of the graph.  However, in the following we explain how the ``average'' prevents this shortcut from working.

\begin{defn}
If $P = (B_{0} \to \ldots \to B_{n})$ is a path in the graph $\Gamma$, then the 
{\em path multiplier}
 is the product of the multipliers along the path, $\lambda(P) := \lambda_{0} \cdots 
\lambda_{{n-1}}$, and the {\em average path multiplier} is $(\lambda(P))^{1/n}$.
\end{defn}

Let $B_k$ and $B_j$ be vertices in a graph $\Gamma$.  Let $P_{k,j}$ be the path from $B_k$ to $B_j$ with the smallest
average path multiplier of all paths from $B_k$ to $B_j$. 
The following example shows that there exist graphs $\Gamma$ such that the cycle containing $B_k$ and $B_j$ with the smallest average multiplier does not contain the path $P_{k,j}$.

\begin{exmp} \label{exmp:avgpathbad}
Again consider the three vertex graph in
Figure~\ref{fig:graphone}, used in the proof of Example~\ref{exmp:avgbad}.  It has two cycles: $0 \to 1 \to 2 \to 0$ and
$0 \to 2 \to 0$, and has multipliers: $\lambda_0 = 1/2,
\lambda_1 = 3/2, \lambda_2 = 8$.  
Comparing paths from vertex $0$ to vertex $2$, we see
$\lambda_0 = 1/2$ and $(\lambda_0 \lambda_1)^{1/2} = \sqrt{3}/2$. The
smallest average path multiplier is along the path $0 \to 2$.  
However,
Example~\ref{exmp:avgbad} shows the smallest
average cycle multiplier is along the cycle $0 \to 1 \to 2 \to 0$. 
\end{exmp}

Example~\ref{exmp:avgpathbad} suggests that 
an algorithm to compute the smallest average cycle multiplier in a graph
must compute the average cycle multiplier
for each simple cycle separately.
A simple combinatorial argument shows that the number of 
simple cycles in a graph (even a sparse graph) is exponential in the number of vertices of the graph.  
The graphs created for polynomial maps of $\CC$ are large enough to make an exponential algorithm prohibitively inefficient (see Section~\ref{sec:examples} for several examples of such graphs).
  Thus, we cannot expect an efficient algorithm for determining the ideal box expansion constant.

The previous exploration can also be used to examine the computational complexity of  Osipenko's method (\cite{Osi00,Osi03}) for approximating Lyapunov exponents, or the Morse spectrum, on $J$, 
and the subsequent hyperbolicity test.    If $C = B_{0} \to \ldots \to B_{{n-1}} \to B_{n}=B_{0}$ is an $n$-cycle of boxes in $\Gamma$, then the characteristic Lyapunov exponent of the cycle is bounded below by 
$$
\frac{1}{n}  \displaystyle\sum_{k=0}^{n-1} \ln ( \lambda_k ) 
= \frac{1}{n} \displaystyle \ln \( \prod_{k=0}^{n-1} \lambda_k \) 
= \frac{1}{n} \ln \( \textit{cycle multiplier of } C \).
$$
Approximating the Morse spectrum of $f$ means computing the characteristic Lyapunov exponent
of each simple cycle in $\Gamma$, and since there are an exponential number of simple cycles,
 this is an exponential-time algorithm.
Osipenko's test for hyperbolicity says that if the minimum exponent over all simple cycles is positive, then the map is hyperbolic.
However, the graph of Example~\ref{exmp:avgpathbad} also suggests there is no shortcut to efficiently computing the minimum exponent, since taking the log and dividing by the cycle length for
this example produces the exact same complications as taking the $n^{th}$ root (the specific numbers can be easily checked). 


\subsection{Approximate, efficient solutions}
\label{sec:approxL}

In order to efficiently determine a successful box expansion constant $L$, we need a good
starting guess.  We have a lower bound of $1$ for $L$.
Corollary~\ref{cor:minmult} gives a way to get an upper bound: if $\mathcal{M}>1$ is the smallest cycle multiplier, realized by a cycle of length $m$, then the ideal box expansion constant (\textit{i.e.}, the smallest average multiplier) $\mathcal{L}$ satisfies $1 < \mathcal{L} \leq \mathcal{M}^{1/m}$.

Alternatively, we can get an upper bound for $L$ in a very simple way.
For a polynomial map $f$, the {\em Lyapunov exponent} 
$\lambda$ measures the rate of growth of tangent vectors to $J$ under 
iteration.  
A description of the one-variable case is given in \cite{Prz}.

\begin{thm} [\cite{Bro, Man}] \label{thm:lyapconnone}
For a polynomial map $f \colon \CC \rightarrow \CC$ of degree $d>1$, the Lyapunov exponent $\lambda$ satisfies $\lambda \geq \log d$,
with equality if and only if $J$ is connected.
\end{thm}

Thus for degree $d$ maps with connected $J$, 
$L=d$ is an 
upper bound.

One straightforward way to obtain a good $L$ in a
preset number of steps is basic \textit{bisection}. Keep track
of lo$L$,  the most recent working L, and hi$L$, the most recent failing $L$.
Lower $L$ halfway to lo$L$ when box expansion fails,
 and raise it halfway to hi$L$ when it succeeds.  
Start with lo$L=1$ and hi$L=d$ or $\mathcal{M}^{1/m}$ if computed.


An alternative to straight bisection is to utilize the
information that Algorithm~\ref{alg:fexpands}, \textbf{BuildMetric}, already discovers in a test for
box expansion.  Indeed, if box expansion fails for some $L$, then we
realized in Proposition~\ref{prop:bestL} that it is due to a cycle, \textit{badcycle},
with
average multiplier $L'$ less than $L$.  But we can easily
adapt the algorithm to compute and return the average
multiplier of \textit{badcycle}, $L'$.  Then if $L' \leq
1$, we know the map is not {\boxexp } on~$\Gamma$ and we can stop.  
Otherwise, on the next pass instead of lowering by some arbitrary amount
we may simply test by this new average multiplier $L'$. Better
yet, we can test by the minimum of~$L'$ and~$L$ minus some
preset step size $\delta$, in order to prevent increasingly small steps
down. 
We refer to the above procedure as the algorithm \textbf{Find-L-Cycles}.

\begin{prop} \label{prop:approxbestL}
Let $\mathcal L$ be the minimum average multiplier over all simple cycles
in the graph $\Gamma$.  If $2 \geq \mathcal L >1+\delta$, then \textbf{Find-L-Cycles}
shows box expansion by
some $L$ within $\delta$ of~$\mathcal L$, in at most
 $1/ \delta$ trials of \textbf{BuildMetric}.
\end{prop}

\begin{proof}
Since we are looking for an expansion amount in
the interval $[1,2]$, and decrease by at least $\delta$ at each step, we
perform at most $1/\delta$ attempts.  Suppose one of these
attempts is successful.  That is, suppose $\Gamma$ fails to box expand
by some $L_0 \leq 2$, and thus outputs an average cycle multiplier
of $L_1 < L_0$. Since $\mathcal L$ is the minimum, $\mathcal L
\leq L_1$. 
Suppose Hypatia verifies successfully box expansion by $L =
\min\{L_{1}, L_0 -\delta \}$.  Then by
Lemma~\ref{lem:hyponeworked} we have $\mathcal L \geq L$. 
Thus we have:
\[  
L_0 > L_1 \geq \mathcal L \geq L=
	\min\{L_{1}, L_0 -\delta \}.
\]
Thus, either we were lucky and $\mathcal L = L_1$ and we have
shown box expansion by exactly that amount, or $L_0 > L_1 >
\mathcal L \geq L_0 - \delta$ and we have shown expansion by
$L_1$ within $\delta$ of~$\mathcal L$.
\end{proof}

The best way to control round-off error in the above with interval arithmetic (see Section~\ref{sec:IA}) is to ``round down'', since a box expansion constant larger than the average cycle multiplier would fail.  That is,  set 
$
L_1 =
\text{Inf}(($Hull$(\lambda_{0}) \cdots $Hull$(\lambda_{{n-1}}))^{1/n}).
$

\section{Examples of running Hypatia for polynomial maps}
\label{sec:examples}

In this section, we describe the results of applying the algorithms of this paper to
 some examples for quadratic and cubic polynomials, $P_c(z)=z^2+c$ and $P_{c,a}(z)
= z^3 -3a^2 z + c$, using our implementation in the computer program Hypatia.
In particular,
in all of the examples of this section, we used the {\namealg } and the method \textbf{Find-L-Cycles} (Section~\ref{sec:approxL}),
to try to produce an $L>1$ and a {\boxmet } for which the constructed $\Gamma$ was \boxexp.
Table~\ref{table:examples} at the end of the section summarizes the data for all 
the examples.

All of the computations described in this section were run on a Sun Enterprise E3500
server with 4 processors, each $400$MHz UltraSPARC (though the 
multiprocessor was not used) and $4$ GB of RAM.  \begin{footnote}{The server was obtained by the Cornell University mathematics department through an NSF SCREMS grant.}\end{footnote}

\subsection{Producing the box model $\Gamma$}

First we summarize how we used the {\boxchcn } from \cite{SLHone} to produce box models $\Gamma$ for these maps.  The {\boxchcn } is in fact an iterative process.  We began by defining a large box $\boxcov_0$ in $\CC$ such that $J \subset \boxcov_0$.  Then for some $n>1$, we place a $2^n \times 2^n$ grid of boxes on $\boxcov_0$.  The construction then builds a graph $\Gamma$ consisting of a subcollection of these grid boxes, which is a box model of $J$, according to Definition~\ref{defn:boxmodel}.  If the boxes are sufficiently small, then every invariant set disjoint from $J$ is disjoint from $\boxcov(\Gamma)$ (for example, attracting periodic orbits).  Further, we can produce an improvement of $\Gamma$ by subdividing all (or some) of the boxes in $\Gamma$, then repeating the construction to produce a new graph $\Gamma'$ such that $\boxcov (\Gamma') \subset \boxcov(\Gamma)$.

For quadratic polynomials, it is easy to check that if $\abs{c} < 2$, then the filled Julia set is strictly contained in the box $[-2,2]^2$.  For the cubic polynomials $P_{c,a}(z)
= z^3 -3a^2 z + c$, one can also check that the filled Julia set is contained in $[-2,2]^2$
whenever $\abs{c} < 2$ and $\abs{a} < (2/3)^{1/2}$.  
We denote $\DD_2 = [-2,2]^2$.
In each of our examples, we began the {\boxchcn } with some $2^n \times 2^n$ grid on $\boxcov_0 = \DD_2$.  

\subsection{Selective subdivision}

We also introduce here two modifications to the basic {\boxchcn}, in order to improve
efficiency. The key idea for both is  \textit{selective subdivision}, \textit{i.e.}, rather than 
subdividing all of the boxes to form the next level, we only subdivide the
boxes where the dynamics is behaving badly. 

Thus we must determine which boxes are most obstructing
the construction of a hyperbolic metric.  One approach is
to choose boxes which are closest to a sink cycle, say 
boxes which have bounded midpoints after several images,
and subdivide only those boxes.  
We call this procedure \textit{\inKselsub}.  See Example~\ref{exmp:cubrabthree}. 

Alternatively, the \textbf{Find-L-Cycles} algorithm,
 in addition to improving the metric on the given
graph, can also be used to determine a selective subdivision
procedure, which we call \textit{\wkcycsub}.  
After \textbf{Find-L-Cycles} is run, we can easily identify the
boxes involved in the cycle with the weakest cycle multiplier.  Then
subdividing these boxes should help improve the expansion amounts.  
We utilize  {\wkcycsub } in Examples~\ref{exmp:bas} and~\ref{exmp:cantcaulmet}.

\subsection{Examples for quadratic polynomials}

The quadratic polynomial $P_c(z) = z^2+c$ has one critical point, at zero, thus has at most one (finite) attracting cycle.


\begin{exmp} \label{exmp:aero}
The aeroplane, $P_c$ for $c=-1.755$, has a period $3$ sink.
Shown in Figure~\ref{fig:hyponeF} 
\begin{figure}
\begin{center}
\drawfighyponeF
\caption{A {\boxJ } for $\boxcov(\Gamma)$, for 
the map $z \mapsto z^2-1.755$, and boxes from a $2^{11} \times 2^{11}$ 
grid on $\DD_2$.  $\Gamma$ is {\boxexp}.}
\label{fig:hyponeF}
\end{center}
\end{figure}
is a {\boxJ } $\boxcov(\Gamma)$
consisting of boxes from $2^{11} \times 2^{11}$ grid on $\DD_2$. This {\boxmod } $\Gamma$
has $43{,}000$ vertices and $260{,}000$ edges.
$\Gamma$ is {\boxexp }, with box expansion constant $L=1.05069$.
The computation took less than $200$ MB of RAM and $6.5$ minutes.
\end{exmp}

\begin{exmp} \label{exmp:bas}
The quadratic polynomial $P(z) = z^2-1$  is called the basilica.
This map has a period $2$ attracting cycle $0 \leftrightarrow -1$.
The picture on the upper left of Figure~\ref{fig:bascant}
is a heuristic sketch of $J$, drawn using the program Fractalasm (available at \cite{CUweb}).   Shading is according to rate of escape to infinity, or to the two-cycle. 

The 
picture on the lower left of Figure~\ref{fig:bascant}
is a {\boxJ } $\boxcov(\Gamma)$, 
composed of selected boxes from
a $2^7 \times 2^7$ grid on $\DD_2$. $\Gamma$ has
$1{,}800$ vertices and $14{,}000$ edges.
We could have easily made a finer picture, but this rough {\boxJ } is enough to prove hyperbolicity. 
Indeed, this {\boxmod } $\Gamma$ is {\boxexp}, with box expansion constant $L=1.14067$.
The associated
{\boxmet } has handicaps, $\varphi_k$, with minimum $0.019$, average $0.049$,
and maximum $1$. 
This initial computation took only $8$ MB of RAM and less
than a minute of CPU time.

For a deeper understanding of this metric, we 
created a picture of the {\boxJ}, with shading of each box according to
the handicap for that box, shown 
 on the lower left of Figure~\ref{fig:bascant}.
Boxes with smaller handicaps are shaded darker.

\end{exmp}

\begin{figure}
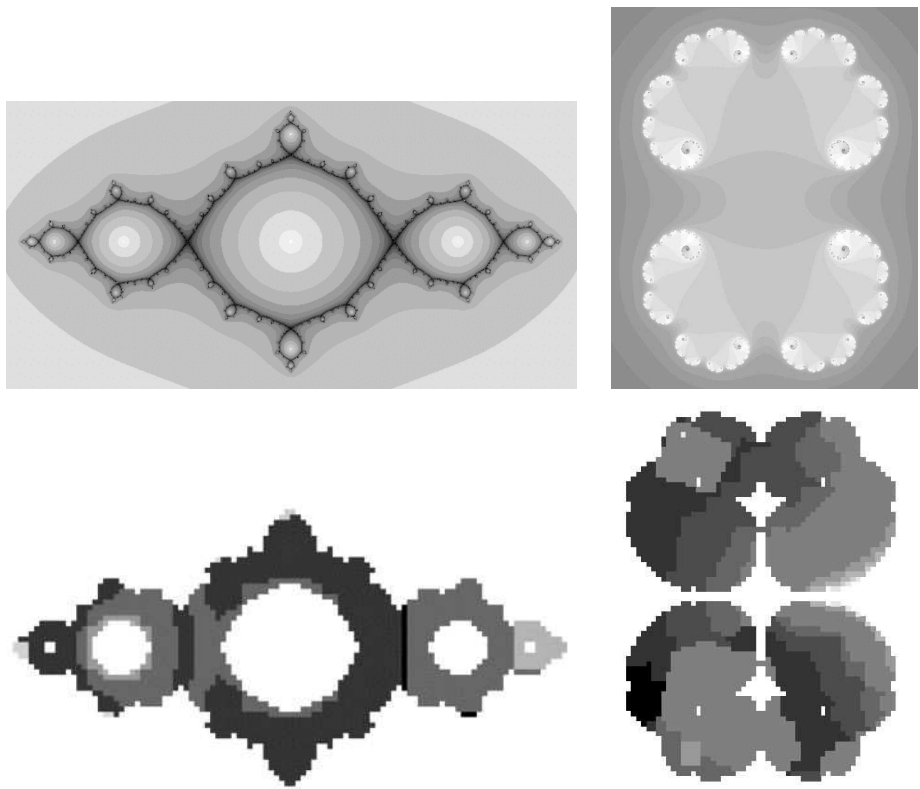

\begin{center}
\drawfigbasFA  \ \ \ \ \drawfigcantcaulFA \\
\drawfigfirstmetex \ \ \ \drawfigcantcaulmet 
\caption
{\textit{Left:} sketches of $J$ for  the map $z \mapsto z^2-1$; 
\textit{Right:} sketches of the cantor $J$ for the map $z \mapsto z^2+.35$;
\textit{Upper:} heuristic sketches of $J$, drawn with the program Fractalasm; 
\textit{Lower:} a box-expansive {\boxJ } $\boxcov (\Gamma)$ for each map, where shading is determined
by handicaps, with smaller handicaps colored darker.
Note the asymmetry in metric shades.  This is not dynamically significant, but an artifact of the tree construction for assigning handicaps.
Boxes on the \textit{lower left} have side length $4/2^7$, and boxes on the \textit{lower right} have side length $4/2^8$.}
\label{fig:bascant}
\end{center}
\end{figure}


\begin{exmp} \label{exmp:cantcaulmet}
A Cantor Julia set near the cusp of the Mandelbrot set is $c=.35$.
The upper right of Figure~\ref{fig:bascant} 
is a Fractalasm sketch of this Julia set.
This is an example in which we can use {\wkcycsub } to show hyperbolicity more quickly.  
The {\boxmod } from the $2^7 \times 2^7$ 
grid on $\DD_2$ fails to be {\boxexp } by $L=1$ due to a bad cycle of length 
only $1$. We had the program subdivide just that box, and again the resulting
{\boxmod } failed for $L=1$ due to a length $1$ bad cycle. 
But, upon subdividing that one box, we found that the resulting {\boxmod } is 
{\boxexp } (by $L=1.00778$).
On the lower right of Figure~\ref{fig:bascant}
is the hyperbolic {\boxJ}, with shading
according to the value of the handicap in each box. Boxes with darker shading have
handicaps closer to zero.
\end{exmp}


\subsection{Examples for cubic polynomials}

The cubics $P_{a,c}(z) = z^3-3a^2z+c$ have two critical points, at $\pm a$, 
thus have at most two attracting cycles.
When $a=0$, there is a clear correspondence between $z^2+c$ and $z^3+c$. 

\begin{exmp} \label{exmp:cubrabthree}
A seemingly rabbit-like cubic is $P_{a,c}, c= -.44-.525i, a=.3i$.
This map has an attracting period three cycle.  
Shown 
in Figure~\ref{fig:cubrabs}
 on the left is the {\boxJ } $\boxcov(\Gamma_{a,c})$ 
with boxes from a
$2^{10} \times 2^{10}$ grid on $\DD_2$.  However, this $\Gamma$ is not hyperbolic because
$\boxcov$ contains one of the critical points, $a=.3i$.  We then
performed {\inKselsub } and produced a {\boxmod } $\Gamma'$ which is {\boxexp}.
The latter {\boxJ } $\boxcov'$ is shown in the center of
Figure~\ref{fig:cubrabs}. 
\end{exmp}

\begin{figure}
\begin{center}
\drawfigcubrabthreea
\drawfigcubrabthreeb
\drawfigcubrabtwo
\caption{Each image above is a box model of $J$ and the period $3$ sink for a cubic polynomial $P_{a,c}$.  Points heuristically found to be in the filled Julia set are shaded lighter,
to illustrate $J$.
\textit{Left:} the {\boxJ } $\boxcov$ for $c= -.44-.525i, a=.3i$, for a 
$2^{10} \times 2^{10}$ grid on $\DD_2$ contains a 
critical point, so the {\boxmod } $\Gamma$ is not {\boxexp}. 
\textit{Center}: a refinement of the left picture, with boxes of side length $4/2^{10}$ and $4/2^{11}$, is {\boxexp}.
\textit{Right:} the {\boxJ } $\boxcov$ for the map $c=-.38125+.40625i, a=0.5i$, with boxes of side length $4/2^9$ and $4/2^{10}$ is {\boxexp}.}
\label{fig:cubrabs}
\end{center}
\end{figure}

\begin{exmp} \label{exmp:cubrabtwo}
The cubic polynomial $P_{a,c},$ with  $c=-.38125+.40625i, a=0.5i$ also has an 
attracting
cycle of period three, but here the Julia set is disconnected.  The most efficient
 method to get a hyperbolic {\boxmod } $\Gamma$ is to
first subdivide all boxes uniformly to obtain a $2^9 \times 2^9$ grid on $\DD_2$, 
and then perform {\inKselsub}.
Shown 
on the right in Figure~\ref{fig:cubrabs}
is the resulting  {\boxJ}. 
\end{exmp}



\begin{sidewaystable}
\caption{Data for the {\boxchmods } described in Section~\ref{sec:examples}.
Here the box depth $n$ is the number such that the boxes are of size $4/2^n$ (from a $(2^n \times 2^n)$ grid on $\DD_2$).}
\label{table:examples}
\begin{center}
{\small
\begin{tabular} { | l | c | c | c | c | c | c | c | c | c | c | c | c | }
\multicolumn{13}{c}{$$}  \\   
\multicolumn{13}{ l }{Examples for the map $P_c(z) = z^2+c$}  \\ \hline\hline

Example &    \multicolumn{2}{c|}{param}   &     sink    &   Figure   &  box   & \# $\Gamma$ & \# $\Gamma$  &    box- & L & max & min & avg  \\ 
	        & \multicolumn{2}{c|}{$c$}    &	 period  &                 &  depth          & boxes    &  edges   
&    exp? &  & $\varphi$   &  $\varphi$  & $\varphi$  \\ 
	        & \multicolumn{2}{c|}{$$}      &	 	     &	                  & 	$n$	        & (1,000s) & (1,000s)  &               &    &        & & \\

\hline\hline

\ref{exmp:aero}    & 
\multicolumn{2}{c|}{$-1.755$}  
& $3$  &   \ref{fig:hyponeF}  &$11$  & $41$ & $250$ & Yes    & $1.0507 $ & $1$ & $2.023\times 10^{-5} $ & $6.8\times 10^{-5}$ 
\\ \hline

\ref{exmp:cantcaulmet} &  
\multicolumn{2}{c|}{$.35$}
 & N/A &  & $7$ & $3.028$  & $27.376$ & No   & $$ & $$ & $$ & $$ 
\\ \cline{5-13}
 &  \multicolumn{2}{c|}{$$}   & & \ref{fig:bascant} R & $7,8,9$ & $3.034$ & $27.422$  & Yes   & $1.00778$ &  $1.19$ &$0.0288723$ & $0.119515$ 
\\ \hline

\ref{exmp:bas}  &   \multicolumn{2}{c|}{$-1$} 
& $2$  &\ref{fig:bascant} L & $7$ &  $1.8$ & $14$ & Yes & 
$1.14067$ & $1$ & $0.019$ & $0.049$ 
\\ \hline

 \multicolumn{13}{c}{$$}  \\ 
\multicolumn{13}{ l }{Examples for the map $P_{a,c}(z)=z^3-3a^2z+c$}  \\ \hline
\hline
Example &    \multicolumn{2}{c|}{params}   &     sink    &   Figure   &  box   & \# $\Gamma'$ & \# $\Gamma'$  &    box- & L & max & min & avg  \\ 
	        &$c$  & $a$  &	 period  &                 &  depth          & boxes    &  edges   
&    exp? &  & $\varphi$  &  $\varphi$   & $\varphi$  \\ 
	        &  &    &	 	     &	                  & 	$n$	        & (1,000s) & (1,000s)  &               &    &        & & \\

\hline\hline

\ref{exmp:cubrabthree} & $-.44$ & $.3i$ & $3$  & \ref{fig:cubrabs} L & $10$ & $92$ & $142$ & No & $$ & $$ & $$ & $$ \\ \cline{5-13}
&$ -.525i$ & & & \ref{fig:cubrabs} C & ${10}, 11$ &  $160$ & $250$ & Yes & $1.0578$ & $1$ & $1.023\times 10^{-4}$ & $3.5\times 10^{-4}$ \\ \hline

\ref{exmp:cubrabtwo}  &   $-.38125$  & $0.5i$ & $3$  & \ref{fig:cubrabs} R   & $9,10$ & $45$ & $867$ &  Yes & $1.0369$ & $1$ & $5.08\times 10^{-4}$  & $.01748$ \\ 
& $+.40625i$ & & & & & & & & & & &  \\  \hline

\hline
\end{tabular}
}  
\end{center}
\end{sidewaystable}

\subsection*{Acknowledgements}
The bulk of this work was accomplished as part of my PhD thesis
at Cornell University (\cite{SLHT}).
I would like to thank John Smillie for providing guidance on this project, 
Clark Robinson for asking a motivating question, 
 James Yorke, Eric Bedford, and John Milnor for advice on the 
preparation of the paper, and Robert Terrell for technical support.


\bibliographystyle{plain}
\bibliography{lynch}

\def\cprime{$'$}
\begin{thebibliography}{10}

\bibitem{CUweb}
Dynamics at~Cornell.
\newblock [http://www.math.cornell.edu/\~{}dynamics].

\bibitem{Bro}
H~Brolin.
\newblock Invariant sets under iteration of rational functions.
\newblock {\em Ark. Mat}, 6:103--144, 1965.

\bibitem{GenIA}
Interval Computations.
\newblock [http://www.cs.utep.edu/interval-comp/].

\bibitem{TCR}
T.~Cormen et~al.
\newblock {\em Introduction to Algorithms}.
\newblock The MIT Electrical Engineering and Computer Science Series. The MIT
  Press and McGraw-Hill Book Company, 1990.

\bibitem{Dell1}
M.~Dellnitz and O.~Junge.
\newblock Set oriented numerical methods for dynamical systems.
\newblock In {\em Handbook of dynamical systems, Vol. 2}, pages 221--264.
  North-Holland, Amsterdam, 2002.

\bibitem{Eiden}
M.~Eidenschink.
\newblock {\em Exploring Global Dynamics: A Numerical Algorithm Based on the
  Conley Index Theory}.
\newblock PhD thesis, Georgia Institute of Technology, 1995.

\bibitem{SLHT}
J.~S.~L. Hruska.
\newblock {\em On the numerical construction of hyperbolic structures for
  complex dynamical systems}.
\newblock PhD thesis, Cornell University, 2002.
\newblock Available for download at the Stony Brook Dynamical Systems Thesis
  Server [http://www.math.sunysb.edu/dynamics/theses/index.html].

\bibitem{SLHthree}
S.~L. Hruska.
\newblock A numerical method for proving hyperbolicity of complex {H}\'{e}non
  mappings.
\newblock submitted, preprint available at [http://xxx.arxiv.org].

\bibitem{SLHone}
S.~L. Hruska.
\newblock Rigorous numerical models for the dynamics of complex {H}\'{e}non
  mappings on their chain recurrent sets.
\newblock submitted, preprint available at [http://xxx.arxiv.org].

\bibitem{Man}
A.~Manning.
\newblock The dimension of the maximal measure for a polynomial map.
\newblock {\em Ann. Math.}, 119:425--430, 1984.

\bibitem{Mil}
J.~Milnor.
\newblock {\em Dynamics in One Complex Variable. Introductory Lectures}.
\newblock Friedr. Vieweg and Sohn, Braunschweig, 1999.

\bibitem{KMisch}
K.~Mischaikow.
\newblock Topological techniques for efficient rigorous computations in
  dynamics.
\newblock {\em Acta Numerica}, 2002.

\bibitem{MooreIA1}
R.~E. Moore.
\newblock {\em Interval Analysis}.
\newblock Prentice-Hall, Englewood Cliffs, New Jersey, 1966.

\bibitem{MooreIA2}
R.~E. Moore.
\newblock {\em Methods and Applications of Interval Analysis}.
\newblock SIAM Studies in Applied Mathematics, Philadelphia, 1979.

\bibitem{Osiold}
G.~Osipenko.
\newblock On the symbolic image of a dynamical system.
\newblock In {\em Boundary value problems (Russian)}, pages 101--105, 198.
  Perm. Politekh. Inst., Perm\cprime, 1983.

\bibitem{Osi}
G.~Osipenko.
\newblock Construction of attractors and filtrations.
\newblock In {\em Conley index theory (Warsaw, 1997)}, volume~47 of {\em Banach
  Center Publ.}, pages 173--192. Polish Acad. Sci., Warsaw, 1999.

\bibitem{Osi00}
G.~Osipenko.
\newblock Spectrum of a dynamical system and applied symbolic dynamics.
\newblock {\em J. Math. Anal. Appl.}, 252(2):587--616, 2000.

\bibitem{Osi03}
G.~Osipenko.
\newblock Symbolic image, hyperbolicity, and structural stability.
\newblock {\em J. Dynam. Differential Equations}, 15(2-3):427--450, 2003.
\newblock Special issue dedicated to Victor A. Pliss on the occasion of his
  70th birthday.

\bibitem{OsiCamp}
G.~Osipenko and S.~Campbell.
\newblock Applied symbolic dynamics: attractors and filtrations.
\newblock {\em Discrete Contin. Dynam. Systems}, 5(1):43--60, 1999.

\bibitem{PBIA}
PROFIL/BIAS Interval~Arithmetic Package.
\newblock \hfil\break
  [http://www.ti3.tu-harburg.de/Software/PROFILEnglisch.html].

\bibitem{Prz}
F.~Przytycki.
\newblock Hausdorff dimension of the harmonic measure on the boundary of an
  attractive basin for a holomorphic map.
\newblock {\em Invent. math.}, 80:161--179, 1985.

\bibitem{War}
W.~Tucker.
\newblock A rigorous {ODE} solver and {S}male's 14th problem.
\newblock {\em Found. Comput. Math.}, 2(1):53--117, 2002.

\end{thebibliography}

\end{document}